\documentclass{elsarticle}

\usepackage{lineno,hyperref}
\journal{Elsevier}

\usepackage[fleqn]{amsmath}
\usepackage{amssymb}

\usepackage{algorithm}
\usepackage{algorithmic}

\usepackage[usenames]{color}

\usepackage{multirow}
\usepackage{array}

\usepackage{moreverb}

\usepackage{graphics}
\usepackage{graphicx}
\usepackage{epsfig}
\usepackage{epstopdf}

\usepackage{footnote}

\usepackage{listings}

\newcommand*\patchAmsMathEnvironmentForLineno[1]{%
  \expandafter\let\csname old#1\expandafter\endcsname\csname #1\endcsname
  \expandafter\let\csname oldend#1\expandafter\endcsname\csname end#1\endcsname
  \renewenvironment{#1}%
     {\linenomath\csname old#1\endcsname}%
     {\csname oldend#1\endcsname\endlinenomath}}% 
\newcommand*\patchBothAmsMathEnvironmentsForLineno[1]{%
  \patchAmsMathEnvironmentForLineno{#1}%
  \patchAmsMathEnvironmentForLineno{#1*}}%
\AtBeginDocument{%
\patchBothAmsMathEnvironmentsForLineno{equation}%
\patchBothAmsMathEnvironmentsForLineno{align}%
\patchBothAmsMathEnvironmentsForLineno{flalign}%
\patchBothAmsMathEnvironmentsForLineno{alignat}%
\patchBothAmsMathEnvironmentsForLineno{gather}%
\patchBothAmsMathEnvironmentsForLineno{multline}%
}

\begin{document}

\begin{frontmatter}

\title{Revisiting Performance of BiCGStab Methods for Solving Systems with
Multiple Right-Hand Sides}

%% Group authors per affiliation:
\author[IMEC]{B.~Krasnopolsky} %\corref{corr_auth}}
%\cortext[corr_auth]{Corresponding author}
\ead{krasnopolsky@imec.msu.ru}

\address[IMEC]{Institute of Mechanics, Lomonosov Moscow State University,
119192 Moscow, Michurinsky ave. 1, Russia}

\begin{abstract}
The paper discusses the efficiency of the classical BiCGStab method and several
of its modifications for solving systems with multiple right-hand side vectors.
These iterative methods are widely used for solving systems with large sparse
matrices. The paper presents execution time analytical model for the time to
solve the systems. The BiCGStab method and several modifications including the
Reordered BiCGStab and Pipelined BiCGStab methods are analyzed and the range of
applicability for each method providing the best execution time is highlighted.
The results of the analytical model are validated by the numerical experiments
and compared with results of other authors. The presented results demonstrate an
increasing role of the vector operations when performing simulations with
multiple right-hand side vectors. The proposed merging of vector operations
allows to reduce the memory traffic and improve performance of the calculations
by about~30\%.
\end{abstract}


\end{frontmatter}

%\linenumbers

%%%%%%%%%%%%%%%%%%%%%%%%%%%%%%%%%%%%%%%

\section{Introduction}
\label{Intro}

The need to solve large sparse systems of linear algebraic equations (SLAEs) is
a typical problem in numerical modelling. The choice of the optimal method,
applicable for the specific problem in terms of compute capacity, memory usage,
scalability in parallel computations and other criteria  can be a challengeable
issue. Nowadays a family of Krylov subspace methods become the popular choice in
a wide range of applications, e.g. the conjugate gradient
(CG,~\cite{Hestens1952}) method for solving SLAEs with symmetric matrices and
the stabilized biconjugate gradient (BiCGStab,~\cite{Vorst1992}) and general
minimal residual (GMRES,~\cite{Saad1986}) methods for solving systems with
general matrices.

The current paper further focuses on the analysis of the
BiCGStab iterative method as an example of the widely used algorithm for solving
SLAEs with matrices of general form. However, the problems discussed in the
paper are common for all Krylov subspace methods, and the proposed ideas can be
extended without in due difficulties to any other iterative method.

\begin{figure}
\centering
\includegraphics[width=5.5cm]{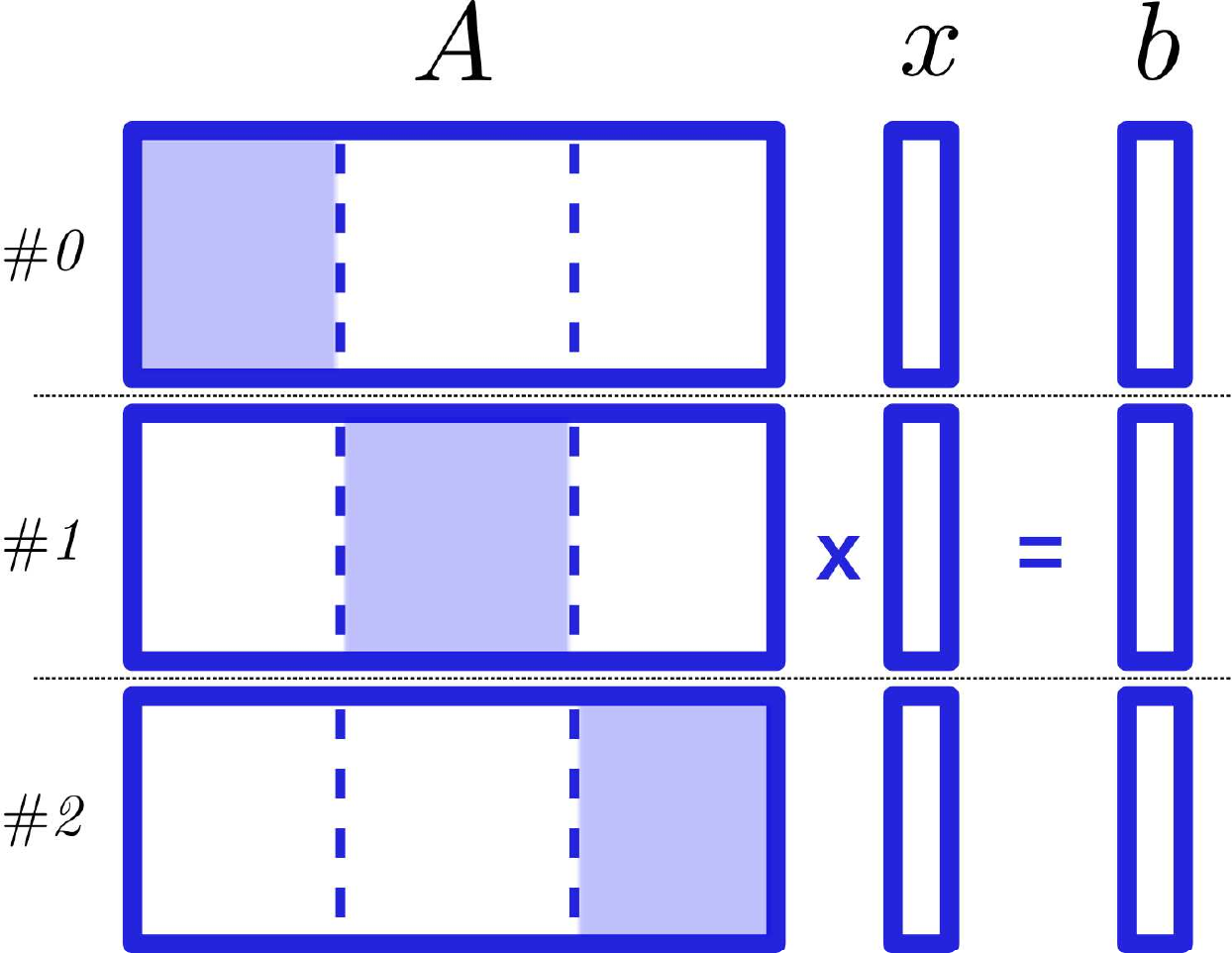}
\caption{Typical data distribution across the computational processes.}
\label{fig:data_distribution}
\end{figure}

The several issues arise when developing efficient parallel implementation of
Krylov subspace iterative methods. Typically the matrix of SLAE and
corresponding vectors are divided row-wise into blocks and distributed across
the computational processes (Figure~\ref{fig:data_distribution}). This allows to
implement vector updates as local operations without any interprocess
communications. The sparse matrix-vector multiplication (SpMV) needs some local
communications between the neighbour processes, but the time for communications
can be overlapped by computations. The third operation, dot product, starts with
the computation of the local vector segments dot products, and continues with
the global reduction among all computational processes. The global reduction leads
to a synchronization point for all computational processes and generally can
affect the overall efficiency of the parallel implementation.

Analyzing the BiCGStab method (Algorithm~\ref{Alg:BiCGStab}
in~\ref{app:merged_methods}) this paper is focused on, one can see that the
method has 2~SpMV operations, 6~vector updates and 4~dot products per each
iteration. In case of preconditioned method the two preconditioning operations
are also performed. Depending on the form of the convergence criteria an
additional dot product may be required to perform the convergence check. The
further discussion assumes the use of criteria in the form of absolute or
relative residual norm, which means the norm of the residual vector $r_{j+1}$
must be additionally computed in line~\ref{BiCGStab_conv_check} of
Algorithm~\ref{Alg:BiCGStab}. Communications for two dot products in
lines~\ref{BiCGStab_dot_merge_1} and~\ref{BiCGStab_dot_merge_2} can be merged,
which reduces the amount of synchronization points per each iteration of
BiCGStab method to~4\footnote{In practice the synchronization point in
line~\ref{BiCGStab_conv_check} can be avoided without any computational overhead
by replacing $(r_{j+1}, r_{j+1}) = (s_j, s_j) - \omega_j (s_j, t_j)$ and
calculating the dot product $(s_j, s_j)$ together with two dot products in
lines~\ref{BiCGStab_dot_merge_1}-\ref{BiCGStab_dot_merge_2} of
Algorithm~\ref{Alg:BiCGStab} in~\ref{app:merged_methods}.}.

The problem of reducing the communications for the Krylov subspace methods, and
specifically the BiCGStab, has received a considerable attention in the
literature. The reduction of communications and the increase of the fraction of
asynchronous data exchanges are achieved by computations reordering and adding
some extra vector operations (vector updates and dot products). This increases
the overall computational costs of the methods but reduces the communication
overhead in parallel simulations. The extra vector operations are usually
ignored when discussing the efficiency of the modified methods as their
computational costs are much lower compared to those of SpMV operations.
This, however, is not the case when the calculations with multiple right-hand
side (RHS) vectors are performed. The input of these operations increases with
increasing the number of RHS vectors, which leads to the need of accurate review
of the corresponding methods.

Most of the modified methods known to date can be classified in the following
categories. The first group of methods focuses on reduction of global
synchronization points and grouping together the dot products. The Modified
BiCGStab method with two global synchronization points has been proposed
in~\cite{MBiCGStab}, and later the Improved BiCGStab method with the single
global synchronization point has been developed~\cite{IBiCGStab}.

The second group of methods focuses on changing the sequence of computations in
order to hide the global communications by the preconditioning
operations or sparse matrix-vector multiplications. The Reordered BiCGStab
method~\cite{RBiCGStab} allows to hide the latency of non-blocking reductions by
performing the preconditioning. This becomes possible due to execution of
preconditioning in advance, i.e. the preconditioning needed at the next
iteration is performed while hiding the global communications at the current
iteration. This leads to an additional single preconditioning operation during
the solution of SLAE. While the proposed modification allows to hide the global
communications, this can provide significant overhead in case of using costly
preconditioner and fast convergence rate of the corresponding system. The
Pipelined BiCGStab method~\cite{Cools2017} extends the same idea and rearranges
the computations in order to allow hiding the global communications behind both
preconditioning and SpMV operation. This makes the method applicable even
without preconditioning, but further increases the computational overhead
compared to the Reordered BiCGStab method, as an additional advanced SpMV
operation is performed. Accounting the different number of extra vector
operations for the Reordered and Pipelined BiCGStab methods it becomes tricky to
determine a priori the area of applicability for each of the methods together
with the classical BiCGStab method.

The third group includes the communication-avoiding BiCGStab
methods~\cite{Carson2013, Naumov2016}, which combine several SpMV operations and
perform the calculation of $s$ iterations an once, thus reducing the amount of
local communications. This approach, however, limits the usage of
preconditioners by the methods of special form, and also adds some extra
operations to reorder the computations.

An additional group of methods, which does not target directly the minimization
of communications problem but discusses important implementation aspects of the
Krylov subspace methods, includes~\cite{Aliaga2013, Anzt2015, Aliaga2015,
Rupp2016}. The authors of these papers focused their research on implementation
of the iterative methods on Graphics Processing Units (GPUs) and demonstrated
the role of kernel fusion technique in the performance boost on GPUs. The
modifications reducing the amount of kernel launch operations for CG, BiCGStab,
and GMRES methods have been proposed. The suggested modifications also reduce
the data transfers between the host and device, which produces an additional
performance boost.

The papers mentioned above were focused on GPU implementation aspects, and
mainly on the single-GPU implementation. The proposed optimizations are valuable
for CPU-based computations, but inapplicable ``as is'' (e.g., the single-GPU
implementation in~\cite{Aliaga2015, Rupp2016} allows to merge the SpMV operation
together with the following vector operations, but in multi-device computations
the data access to the resulting vector is indirect and can not be efficiently
merged with subsequent vector operations). The modifications of the methods
proposed in~\cite{Anzt2015} presume the more robust formulation and preserve the
sparse matrix-vector multiplication as an independent function call. This
decreases the overall performance gain for the loop fusion technique compared
the fully fused formulations, but ensures the algorithm does not depend on the
matrix storage format and the SpMV operation implementation details.

The current paper attempts to systematically revise the BiCGStab method and
several modified versions including the Improved BiCGStab, Reordered BiCGStab
and Pipelined BiCGStab in terms of execution time and parallel efficiency with
the focus on the computations with multiple right-hand side vectors (the
possible convergence issues due to reordering the operations are not accounted
in this paper). The rest of the paper is organized as follows. The impact of
reducing the amount of data transfers with the memory and merging the vector
operations on the performance of BiCGStab methods is discussed in the second
section. The third section is devoted to construction of data transfer-based
execution time model. The model is applied to the considered BiCGStab methods
and validated by the corresponding simulation results. The detailed comparison
of the BiCGStab methods is presented in the fourth section. The specifications
of compute platforms used for the calculations and merged formulations of the
BiCGStab methods are summarized in Appendix.

%%%%%%%%%%%%%%%%%%%%%%%%%%%%%%%%%%%%%%%%%%%%%%%%%%%%%%%%%%%%%%%%%%%%%%%%%%%%%%%

\section{Merging vector operations in BiCGStab methods}
\label{MergedSolvers}

\subsection{The role of merging vector operations on CPUs}
The loop fusion technique applied for GPUs was focused on reduction of data
transfers between the host and device memory and decrease of the time losses due
to multiple kernel launches. This technique, however, is also of importance for
the CPU computations. Merging of several vector operations in some cases allows
to reuse the data already transferred from the memory and reduce the total
volume of data transfers. Taking into account that the performance of basic
operations comprising the iterative methods for solving SLAEs with sparse
matrices is limited by the memory bandwidth of the compute
system~\cite{Williams2009}, this allows to increase the compute intensity of the
algorithm (flop per byte ratio) and, consequently, performance of the
calculations.

An importance of the corresponding modifications can be illustrated by the
simple example performing the calculation of vector updates: $y = a x + b y$, $z
= b x + a z$, and $z = b y + c z$, where $a$, $b$, and $c$~are the scalars and
$x$, $y$, and $z$~are the vectors of size~$N$. The corresponding operations are
computed in three ways: (1)~three independent BLAS-like function calls,
(2)~merged loop performing two vector updates, and (3)~merged loop performing
three vector updates (Figure~\ref{fig:merge}). The first run needs $9N$~elements
to be transferred with the memory and $9N$~floating point operations (FLOP), the
second run transfers $5N$~elements and performs $6N$~FLOP, and the last one
transfers $5N$~elements with $9N$~FLOP. The corresponding calculation times for
the double precision floating point vectors of total length $N = 10^8$,
performed on Lomonosov and Lomonosov-2 compute systems, are shown in
Tab.~\ref{tab:axpby} (characteristics of the test platforms can be found in
Tab.~\ref{tab:spec} of~\ref{section:spec}). The presented results are obtained
for the single nodes and use MPI with ``one rank per core'' scheme to utilize
all processor cores. One can see the ratio of the execution times to perform the
first and second cases is about~1.76. This value is closer to the memory traffic
ratio~1.8 than the FLOP ratio~1.5. Moreover, the third case produces the same
memory traffic as the second one, but performs 1.5~times as much floating point
operations, however, the execution time equals to the second case. This clearly
demonstrates that the execution time for these computations depends on the
memory traffic but not on the number of floating point operations, and the
overhead to compute the third vector update for the data already loaded to the
registers is negligible.

\begin{figure}[t!]
\begin{lstlisting}
// run #1:
for(i = 0; i < N; i++)
  y[i] = a * x[i] + b * y[i];
for(i = 0; i < N; i++)
  z[i] = b * x[i] + a * z[i];
for(i = 0; i < N; i++)
  z[i] = b * y[i] + c * z[i];

// run #2:
for(i = 0; i < N; i++) {
  y[i] = a * x[i] + b * y[i];
  z[i] = b * x[i] + a * z[i];
}

// run #3:
for(i = 0; i < N; i++) {
  y[i] = a * x[i] + b * y[i];
  z[i] = b * x[i] + a * z[i];
  z[i] = b * y[i] + c * z[i];
}
\end{lstlisting}
\vspace{-0.5cm}
\caption{Pseudocode of the three groups of vector operations used to analyze
the impact of memory traffic on performance of the vector operations.}
\label{fig:merge}
\end{figure}

\begin{table}[t!]
  \caption{Time to perform vector operations, ms.} \label{tab:axpby}
\centering
\begin{tabular}{ | c | c | c | }
\hline
	 & Lomonosov & Lomonosov-2 \\
\hline
	run \#1 & 207 & 125 \\
\hline
	run \#2 & 117 & 71 \\
\hline
	run \#3 & 116 & 71 \\
\hline
\end{tabular}
\end{table}

%%%%%%%%%%%%%%%%%%%%%%%%%%%%%

An observation demonstrated above suggests two conclusions: (1)~it is reasonable
to reformulate the corresponding iterative methods with the focus on the
minimization of memory transfers and grouping vector operations, and (2)~the
execution time model must include the metrics based on the volume of data
transfers, but not the floating point operations.

%%%%%%%%%%%%%%%%%%%%%%%%%%%%%

\subsection{Merged formulations of BiCGStab iterative methods}

The simple example presented in the previous paragraph clearly shows importance
of merging the vector operations in Krylov subspace methods to achieve maximal
performance. The current paragraph provides the merged formulations for the
modified BiCGStab methods discussed in the rest of the paper.
It should be noted the applicability of the approach considered
is not limited by the chosen group of methods and can also be applied to the
other iterative methods.

The following principles make the basis of constructing the merged formulations
of the methods. The merging of vector operations is applied if:
\begin{itemize}
\item merging does not increase the amount of communications and overall
computational costs of the method (i.e., no extra computations allowed after the convergence
criteria satisfied);
\item merging does not decrease the amount of computations that can be used to
hide the asynchronous communications;
\item merged operation reduces the number of vector reads/writes.
\end{itemize}
In contrast with GPU computations, CPUs do not have observable penalty on
additional function calls, and reduction of overall function calls has
negligible impact on the overall performance.

Considering the BiCGStab method presented in Algorithm~\ref{Alg:BiCGStab}
in~\ref{app:merged_methods}, three groups of merged vector operations are
suggested. Grouping of three dot products in line~\ref{MBiCGStab_group1} of
Algorithm~\ref{Alg:MBiCGStab} in~\ref{app:merged_methods}, needed to compute the
coefficient $\omega_j$ and norm of the residual vector $r_{j+1}$, allows to
reduce the memory traffic from six vector reads to only two. Grouping of two
vector updates in line~\ref{MBiCGStab_group2} allows to reduce single vector
read operation (this branch, however, is executed only once and has little
effect on the overall performance). Finally, grouping in
line~\ref{MBiCGStab_group3} allows to reduce two vector reads. In total, the
memory traffic for the single loop of the BiCGStab method decreases from 18~read
and 4~write operations to 14~reads and 4~writes in the merged formulation.

Applying the same principles to the modified BiCGStab methods one can obtain the
merged formulations for these methods. The corresponding algorithms are
presented in~\ref{app:merged_methods} and the data traffic characteristics of
the proposed formulations are summarized in Tab.~\ref{tab:loop_mem_traffic}.
While merging of vector operations for classical BiCGStab method allows to
reduce data traffic for vector operations by only 18\% and this value is even
lower for the preconditioned method, the effect of merging is much stronger for
the modified formulations: for Improved BiCGStab and Pipelined BiCGStab methods
the data traffic can be decreased twofold. Merging of vector operations allows to
significantly reduce the overhead occurred as a result of reordering the
computations to decrease the communication costs.

\begin{savenotes}
\begin{table}[t!]
  \caption{Number of vector read/write operations per single iteration of the
  iterative method.\label{tab:loop_mem_traffic}}
\centering
 \begin{tabular}{ | c | c | c | c | c | c | c | c |}
 \hline
  \multirow{2}{*}{Method} & \multicolumn{3}{|c|}{Basic} & \multicolumn{3}{|c|}{Merged} \\
  \cline{2-7}
         & Read & Write & Total & Read & Write & Total \\
\hline
	BiCGStab & 18 & 4 & 22 & 14 & 4 & 18 \\
\hline
 	IBiCGStab & 27 & 6 & 33 & 14 & 6 & 20 \\
\hline
 	PipeBiCGStab\footnote{Formally, the Pipelined BiCGStab method needs 34~reads and
8~writes. However, the vector operation in line~\ref{PipeBiCGStab_4vec} of
Algorithm~\ref{Alg:PipeBiCGStab} assumes the vector update including 4~vectors.
The standard BLAS functions perform operations with only two vectors, and the
extended BLAS functions allow three vector arguments~\cite{BLAS_standard}.
Without writing the specific BLAS-like function with 4~vectors this vector
update can be implemented with at least two function calls leading to additional
vector read and write operations.}
 	 & 34 & 9 & 43 & 18 & 8 & 26 \\
\hline
	PBiCGStab & 18 & 4 & 22 & 15 & 4 & 19 \\
\hline
	RBiCGStab & 22 & 6 & 28 & 15 & 6 & 21 \\
\hline
	PPipeBiCGStab\footnote{The preconditioned Pipelined BiCGStab method has two vector
operations with 4~vectors.} 
	 & 43 & 13 & 56 & 23 & 11 & 34 \\
\hline
\end{tabular}
\end{table}
\end{savenotes}

%%%%%%%%%%%%%%%%%%%%%%%%%%%%%%%%%%%%%%%

\subsection{Numerical experiments}
\label{sec:numerical_merge}

The presented theoretical prepositions are validated by the numerical simulation
results. The corresponding algorithms are implemented in a newly developing XAMG
library of numerical methods for solving SLAEs with multiple RHS vectors. The
testing is done for the matrix obtained as a result of spatial discretization of
the Poisson equation in a cubic domain using 7-point stencil with the uniform
grid of $200^3$ cells and single RHS vector (of 8~mln. unknowns). The fixed
number of iterations is used for the benchmark purposes, $N_{it} = 1000$, to
measure the corresponding execution time. The simulations are performed on the
single compute node and utilize all available processor cores; the computations
are parallelized with help of MPI, mapping one MPI rank per each processor core.
In case of preconditioned methods the identity operator is used as a
preconditioner, which leads to a simple vector copy, equivalent in terms of data
traffic to a single vector read and write operations (the preconditioning costs
are not included in Tab.~\ref{tab:loop_mem_traffic}).

The testing is performed for six iterative methods including classical BiCGStab,
Improved BiCGStab, Pipelined BiCGStab, preconditioned BiCGStab, Reordered
BiCGStab and preconditioned Pipelined BiCGStab methods for both basic and merged
formulations. The corresponding execution times are summarized in
Tab.~\ref{tab:loop_exec_time}. The classical BiCGStab method demonstrates the
lowest execution time. The Improved BiCGStab takes 20\% more time for basic
formulation and 7\% for the merged one. The Pipelined BiCGStab has overhead of
about 40\%, but thanks to merging of vector operations it can be reduced
twofold.

The preconditioned BiCGStab is also the fastest among the preconditioned
methods. The Reordered BiCGStab shows 7\% overhead, and it becomes as small as
2\% with the merged formulation. The preconditioned Pipelined BiCGStab slows
down by 60\%, but this value can also be reduced twofold if using merged
formulation of the method.

The presented results demonstrate that merging of vector operations is a useful
optimization for CPU computations allowing to decrease the calculation time by
reducing the data traffic. The execution time reduction varies depending on the
method from 4\% to 21\%. The proposed optimization decreases the computational
costs to perform extra vector operations introduced as a result of reordering
thus allowing to extend the range of applicability of modified methods.

The efficiency of the implemented methods is compared with the open-source
\textit{hypre}~\cite{hypre} library (\textit{hypre} v.2.15.1 has been built with
Intel compiler and Intel MKL, and with MPI and OpenMP parallelization). The
measured execution times for the unpreconditioned BiCGStab method are up to 20\%
higher than the ones obtained with the current implementation of basic BiCGStab
method formulation. Compared the merged formulations the overall speedup reaches
36\% and 28\% for Lomonosov and Lomonosov-2 supercomputers correspondingly.

\begin{table}[t!]
  \caption{Execution times for the single iteration
  of the numerical method, times in ms.\label{tab:loop_exec_time}}
\centering
 \begin{tabular}{ | c | c | c | c | c | c | c | c |}
 \hline
  \multirow{2}{*}{Method} & \multicolumn{2}{|c|}{Lomonosov} & \multicolumn{2}{|c|}{Lomonosov-2} \\
  \cline{2-5}
         & Basic & Merged & Basic & Merged \\
\hline
	BiCGStab & 113 & 105 & 61.5 & 57.5 \\
\hline
 	IBiCGStab & 136 & 113 & 74.3 & 61.2 \\
\hline
 	PipeBiCGStab & 157 & 126 & 86.5 & 68.6 \\
\hline
	PBiCGStab & 121 & 116 & 66.4 & 63.7 \\
\hline
	RBiCGStab & 129 & 117 & 72.1 & 64.8 \\
\hline
	PPipeBiCGStab & 195 & 155 & 107 & 84.8 \\
\hline
	\textit{hypre}, BiCGStab (MPI) & 143 & -- & 81.8 & -- \\
\hline
	\textit{hypre}, BiCGStab (MPI+OpenMP) & 144 & -- & 74 & -- \\
\hline
\end{tabular}
\end{table}

The effect of merging vector operations for the simulations with multiple
RHS~vectors is demonstrated in Tab.~\ref{tab:loop_exec_time_nv}. Using the same
test matrix, the series of experiments are repeated for SLAEs with 4 and 16 RHS
vectors. The obtained results show the increase of the merging effect when
performing simulations with multiple RHS vectors, and the higher a method is
affected by vector merging, the stronger an effect appears to be. While the
classical BiCGStab shows increase by only~3\%, the use of vector merging for the
Pipelined BiCGStab allows to accelerate the simulations with 16~RHS vectors by
27\%.

\begin{table}[t!]
  \caption{Execution times and corresponding improvement due to merging of
  vector operations for the single iteration of the numerical method with
  variable number of RHS vectors. Lomonosov-2 supercomputer, times in seconds.
  \label{tab:loop_exec_time_nv}}
\centering
 \begin{tabular}{ | c | c | c | c | c | c | c | c | c | c |}
 \hline
  \multirow{2}{*}{Method} & \multicolumn{3}{|c|}{$m = 1$}
  & \multicolumn{3}{|c|}{$m = 4$} & \multicolumn{3}{|c|}{$m = 16$} \\
  \cline{2-10}
         & Basic & Merged & Imp. & Basic & Merged & Imp. & Basic & Merged &
         Imp.\\
\hline
	BiCGStab & 0.062 & 0.058 & 7\% & 0.181 & 0.164 & 9\% & 0.65 & 0.59 & 10\%\\
\hline
 	IBiCGStab & 0.074 & 0.061 & 18\% & 0.236 & 0.182 & 23\% & 0.87 & 0.66 & 25\%\\
\hline
 	PipeBiCGStab & 0.087 & 0.069 & 21\% & 0.282 & 0.212 & 25\% & 1.06 &
 	0.78 & 27\% \\
\hline
	PBiCGStab & 0.066 & 0.064 & 3\% & 0.201 & 0.189 & 6\% & 0.72 & 0.68 & 7\%\\
\hline
 	RBiCGStab & 0.072 & 0.065 & 10\% & 0.224 & 0.196 & 12\% & 0.81 & 0.7 & 14\%\\
\hline
 	PPipeBiCGStab & 0.107 & 0.085 & 21\% & 0.368 & 0.279 & 24\% & 1.39 & 1.03 &
 	26\% \\
\hline
\end{tabular}
\end{table}

%%%%%%%%%%%%%%%%%%%%%%%%%%%%%%%%%%%%%%%%%%%%%%%%%%%%%%%%%%%%%%%%%%%%%%%%%%%%%%%

\section{Execution time model}
\label{Model}

A variety of modifications of numerical methods published in the literature
leads to the need of formulating simplified models to classify them and
highlight the range of applicability for each of the methods. The corresponding
attempts to construct analytical models have been done earlier, e.g.,
in~\cite{deSturler1994, Zhu2014}. The basic execution time model has been
proposed in~\cite{deSturler1994}. The model includes the time to perform the
floating point operations and the time of global communications for dot
products, but ignores the time spent on SpMV communications. Using this model
the authors proposed parallel modifications of CG and GMRES methods, focused on
reduction of global communications and hiding communications by local
computations.

The corresponding model has been revised and applied to Generalized Product type
of Bi-Conjugate Gradient method~\cite{ref:Zhang1997, Fujino2002}
in~\cite{Zhu2014}. This model is also based on the amount of floating point
operations and accounts the time spent on data exchanges in SpMV operations and
dot product calculations. The corresponding communication time dependencies on
the number of processes are expressed as a power functions, and the free
parameters are calibrated with results of numerical experiments.

The major drawback of these models relies in the usage of floating point
operations as a metric of the execution time. As shown in the previous section,
the basic operations in Krylov subspace methods applied to solve large sparse
SLAEs are memory bound and limited by the memory bandwidth, but not by the
compute capacity of modern processors. Thus, the execution time for the
operations is determined by the amount of data to be read from and written to
the memory, but not by the corresponding floating point operations.

%%%%%%%%%%%%%%%%%%%%%%%%%%%%%

\subsection{Data transfer-based execution time model}

The data transfer-based execution time model is constructed to analyze the
efficiency and range of applicability for various BiCGStab methods.
The model covers simulations with multiple RHS vectors and additionally allows
to account various preconditioners. It is assumed multiple vectors are stored
row-wise in a single vector, allowing to improve the cache data reuse when
performing SpMV operations and utilize the vectorization effects.
The minimal compute unit in the analytical model is chosen equal to the single
compute node.

The execution time for the single iteration loop of the BiCGStab methods
includes the times to perform local vector operations, $T_{vec}$, global data
exchanges when calculating dot products, $T_G$, sparse matrix-vector
multiplications, $T_{mul}$, and communications with local neighbours when
computing SpMVs, $T_L$. In case of preconditioned methods, an additional time
spent on preconditioning, $T_{prec}$, must be included. The time for local
calculations can be expressed as a ratio of the volume of data transfers,
$\Sigma$, to the compute system memory bandwidth, $B$:
\begin{gather}
T = \frac{\Sigma}{B}, \label{eq:exec_time}
\end{gather}
where the overall memory bandwidth is a multiplication of the memory bandwidth
of the single node, $b$, to the number of compute nodes, $p$, used in the
calculations:
\begin{gather}
B = b \, p.
\end{gather}

Assuming the integers and floating point numbers are of size 4 and 8~bytes
respectively, and the sparse matrix is stored in CRS format~\cite{Saad2003}, the
corresponding sparse matrix-vector multiplication
needs~\cite{KrasnopolskyCPC2018}
\begin{gather}
\Sigma_{mul} = N \left( 8m (C+1) + 4 (3C+1) \right) \label{eq:SpMV_mem}
\end{gather}
bytes to be transferred with the memory. Here, $N$ is the matrix size, $m$ is
the number of RHS vectors, and $C$ denotes the average number of nonzeros per
matrix row. This leads to the SpMV local computations execution time
\begin{gather}
T_{mul}(p) = \frac{N \left(8 m (C+1) + 4(3C+1) \right)}{b \, p}.
\end{gather}

The memory traffic of single vector read/write operation is expressed as
\begin{gather}
\Sigma_{vec} = 8 \, N \, m.
\end{gather}
Substituting this expression to~\eqref{eq:exec_time}, the estimate for the
vector operations execution time takes the form:
\begin{gather}
T_{vec}(p) = \frac{8 N m}{b \, p}.
\end{gather}
To complete, the local and global communication times as a functions of the
number of compute nodes and message size must be provided. 

%%%%%%%%%%%%%%%%%%%%%%%%%%%%%%%%%%%%%%%

\subsubsection{Local communications}

The local non-blocking point-to-point communications arise when performing data
exchange to complete the SpMV operation. Generally, the communication pattern
(list of neighbour processes and message sizes) depends on a couple of factors
including the matrix topology and data decomposition, and construction of the
corresponding time estimates can be a challenging issue. In order to simplify
the further narration and focus on the effect of global communications the
specific form of the test matrices is considered for parallel runs. The case
assumes the communication pattern remains close to constant with increasing the
number of compute processes. Such a scenario is realized, for example, for
banded matrices with narrow band.

The corresponding communication pattern can be modelled with the simple
benchmark, performing asynchronous data exchange for each process with only
several neighbours. This test is a generalization of IMB Exchange
benchmark~\cite{IntelMPIBench}, which allows to measure the latency of the
non-blocking point-to-point communications depending on the number of neighbours
and transferred message size. The testing is performed on Lomonosov
supercomputer using up to 128~compute nodes and 8~MPI processes per node.
The obtained results show that the time to transfer the data is almost
independent on the overall number of MPI processes and weakly depends on the
number of neighbours if it does not exceed~10. The local communications time can
be expressed as a function of the transferred message size in the form
\begin{gather}
T_L(l) = C_0 + C_1 \, l^{n_0},
\end{gather}
where $l$ is the message size in bytes. The corresponding data fitting of the
benchmark results gives two sets of coefficients depending on the message size:
\begin{gather}
T_L(l) = \begin{cases} 
2.4 \cdot 10^{-6} + 6.9 \cdot 10^{-8} \, l^{0.56} &, \, l \leq 2048 \,\,
\mbox{bytes}, \\
3.2 \cdot 10^{-6} + 2 \cdot 10^{-9} \, l &, \, l > 2048 \,\, \mbox{bytes}.
\end{cases}
\end{gather}

%%%%%%%%%%%%%%%%%%%%%%%%%%%%%%%%%%%%%%%

\subsubsection{Global communications}
\label{global_comm}

The need for global communications arise after performing the dot products with
local vector fragments. The corresponding operation using the MPI library to
perform communications can be implemented in the form of blocking call of
\textit{MPI\_Allreduce} or non-blocking call of \textit{MPI\_Iallreduce}
function. Theoretically, the second one can be preferable for the algorithms
with overlap of global communications by computations (e.g., the Reordered
BiCGStab or Pipelined BiCGStab methods). In practice, however, the situation is
not so evident. The efficiency of the non-blocking global communications depends
on lots of factors, including the MPI library implementation and communications
hardware.
To perform the non-blocking collective operations (e.g.,
\textit{MPI\_Iallreduce}) the specific progression must be applied. The MPI
standard~\cite{MPI_standard} does not specify the progression rule, and its
implementation varies depending on the specific MPI
library~\cite{Cardellini2016}. Three basic progression strategies are discussed
in the literature~\cite{Hoefler2008}: hardware-based progression, software-based
progression, and manual progression. The first option seems to be the most
promising one, but is available in a limited number of compute systems (e.g.,
Cray systems with Gemini or Aries interconnect~\cite{Cardellini2016}), and not
fully supported in the majority of other systems. The software based progression
is implementation-specific, and not limited by the communication hardware. This
option is typically implemented with threads, which handle the status of the
non-blocking operations and perform the corresponding progression. The drawback
for this strategy is related with significant overhead, produced by the
progression threads~\cite{Hoefler2008, Si2018, Ruhela2018, Denis2019}. The
manual progression is generally independent on the hardware and MPI library
implementation, but needs some user efforts to add \textit{MPI\_Test} or
\textit{MPI\_Probe} calls to progress the communications. Each of the strategies
mentioned above has its advantages and disadvantages depending on the performed
communication and message size. The corresponding discussion can be found in,
e.g.~\cite{Cardellini2016, Hoefler2008, Si2018, Ruhela2018, Denis2019}.

The efficiency of global communications and various progression techniques for
the Lomonosov supercomputer has been investigated in detail with help of
corresponding benchmark in~\ref{app:Iallreduce_async}. The
efficiency of overlap of communications and computations is estimated in terms
of overlapping overhead, $\gamma$. This parameter is defined as
\begin{gather}
\gamma = \frac{T_{calc}^{async} + T_{comm}^{async} -
T_{calc}^{sync}}{T_{comm}^{sync}},
\end{gather}
where $T_{calc}$ is the calculation time and $T_{comm}$ is the time spent on
communications. Superscripts ``sync'' and ``async'' correspond to blocking
synchronous and non-blocking asynchronous communications respectively;
$T_{calc}^{async}$ can be higher than $T_{calc}^{sync}$ due to influence of
communications performing in background with calculations (the details can be
found in~\ref{app:Iallreduce_async}). The parameter introduced characterizes the
degree of overlap of communications by computations with account of potential
increase of the calculation time due to progress of asynchronous non-blocking
communications.

The presented results demonstrate that the special efforts are needed to obtain
the overlap of non-blocking global communications by computations. The
experiments with Intel MPI~2017 library have shown that, while the software
progression provides almost ideal overlap with overlapping overhead of only
10-20\%, the execution times for the global reductions with short messages (less
than 8192~bytes) increase more than by the order of magnitude compared the ones
without software progression. This makes the corresponding functionality
inapplicable in real simulations. The manual progression allows to obtain
overlapping overhead by about 30-50\%, however, it requires to perform
systematic \textit{MPI\_Test} progression calls and it can be tricky to
implement this progression in real computational codes.

The simple analytical expression for blocking global communications can be
constructed as a result of generalization of benchmark results presented
in~\ref{app:Iallreduce_async}. Following~\cite{Zhu2014}, the power function is
used to approximate the measured data. It is found the measured data can be
approximated with acceptable degree of accuracy by the function
\begin{gather}
T_{G}(p, l) = C_0 + C_1 \, l^{n_0} \, p^{n_1}, \label{eq:global_comm}
\end{gather}
where the message size $l$ is equal to $l = 8 \, m \, k$, and $k$ is the number
of dot products computed at once. For the iterative methods considered in the
paper this parameter varies from~1 to~7. The corresponding coefficients $C_0$,
$C_1$, $n_0$, and $n_1$ are defined by fitting results of numerical
experiments. For Lomonosov supercomputer these coefficients are equal to:
\begin{gather}
C_0 = 3.5 \cdot 10^{-6}, \, C_1 = 1.7 \cdot 10^{-6}, \, n_0 = 0.21, \, n_1 =
0.54.
\end{gather}

The overlapping overhead parameter, $\gamma$, is used to account the degree of
overlap when performing asynchronous global communications. The corresponding
execution time for the calculations and non-blocking asynchronous communications
takes the form:
\begin{gather}
T = \max(T_{calc} + \gamma T_G, T_G). \label{eq:overlap}
\end{gather}

%%%%%%%%%%%%%%%%%%%%%%%%%%%%%%%%%%%%%%%

\subsubsection{Execution times for the BiCGStab methods}

Using the expressions introduced above one can obtain the expected execution
times for the iterative methods. The corresponding relations presented below
operate with merged formulations of the methods with reduced number of vector
reads/writes. The classical BiCGStab method performs all the calculations
step-by-step, except the calculation of sparse matrix-vector multiplication and
local communications, which can be efficiently overlapped. This leads to the
following form of the theoretical execution time for the single iteration of the
BiCGStab method:
\begin{equation}
T^{BiCGStab}(p) = 18 T_{vec}(p) + 2 T_{SpMV}(p) + T_G(p, 24m) + \\
2T_G(p, 8m), \label{est:bicgstab}
\end{equation}
where
\begin{equation}
T_{SpMV}(p) = \max\left( T_{mul}(p), T_L(l) \right) \label{est:spmv}
\end{equation}
and $l$ is the average message size transferred with the neighbour processes.
Opposite to the global communcations, it is assumed that local non-blocking
asynchronous point-to-point communications do not produce observable overlapping
overhead~\cite{Medvedev2019} and this overhead is ignored in the analytical
model. The real average message size depends on the specific test matrix and the
corresponding data distribution across the compute processes. For the sake of
simplicity, in the following discussion the message size is assumed independent
on the number of compute processes utilized. This assumption is valid for the
test matrix used in the following validation section. Meanwhile, the
corresponding dependence can be included without in due difficulties to the
proposed execution time model.

The Improved BiCGStab method generally has the same structure of the algorithm
as the classical method formulation with the only difference in the number of
vector operations and global communications:
\begin{equation}
T^{IBiCGStab}(p) = 20 T_{vec}(p) + 2 T_{SpMV}(p) + T_G(p, 56m).
\label{est:ibicgstab}
\end{equation}
The message size for the global communications in the Improved BiCGStab is about
twice higher than the ones in the classical BiCGStab method, but changed order
of computations allows to decrease the number of global reductions, which is
typically preferable in the high-scale simulations. The Pipelined BiCGStab
method allows to hide the global communications by the SpMV operation.
Accounting the overlap of global communications and computations, the following
relation for the execution time can be formulated:
\begin{multline}
T^{PipeBiCGStab}(p) = 26 T_{vec}(p) + \\
\max \left( T_{SpMV}(p) + \gamma T_G(p, 24m), T_G(p, 24m)\right) + \\
\max \left( T_{SpMV}(p) + \gamma T_G(p, 32m), T_G(p, 32m)\right).
\label{est:pipebicgstab}
\end{multline} 

The preconditioned methods reproduce the same sequence of computations and
additionally perform two preconditioning operations per each iteration of the
methods. For the preconditioned BiCGStab method the corresponding times are not
overlapped with any computations and simply added to the overall execution time:
\begin{multline}
T^{PBiCGStab}(p) = 19 T_{vec}(p) + 2 T_{SpMV}(p) + 2 T_{prec}(p) + T_G(p, 24m) + \\
2 T_G(p, 8m).\label{est:pbicgstab}
\end{multline}
The Reordered BiCGStab method allows to overlap preconditioning operations with
global reductions. Using expression~\eqref{eq:overlap} the corresponding
expression takes the form:
\begin{multline}
T^{RBiCGStab}(p) = 21 T_{vec}(p) + 2 T_{SpMV}(p) + \\
\max \left(T_{prec}(p) + \gamma T_G(p, 8m), T_G(p, 8m)\right) +  \\ 
\max \left(T_{prec}(p) + \gamma T_G(p, 32m), T_G(p, 32m)\right).
\label{est:rbicgstab}
\end{multline}
Finally, the preconditioned Pipelined BiCGStab method allows to overlap global
communications with preconditioning and sparse matrix-vector multiplication:
\begin{multline}
T^{PPipeBiCGStab}(p) = 34 T_{vec}(p) + \\ 
\max \left( T_{SpMV}(p) + T_{prec}(p) + \gamma T_G(p, 24m), T_G(p, 24m)\right) + \\ 
\max \left( T_{SpMV}(p) + T_{prec}(p) + \gamma T_G(p, 32m), T_G(p, 32m)\right).
\label{est:ppipebicgstab}
\end{multline}

%%%%%%%%%%%%%%%%%%%%%%%%%%%%%%%%%%%%%%%%%%%%%%%%%%%%%%%%%%%%

\subsection{Execution time model validation}

The proposed execution time model is thoroughly validated and compared with the
corresponding calculation results. Validation includes the single-node runs,
performed on two compute systems, and multi-node runs, performed on 128~compute
nodes of Lomonosov supercomputer.

\subsubsection{Single-node validation}

The single-node validation is performed for the test matrix of 8~mln. unknowns,
previously used in section~\ref{sec:numerical_merge}. This test series includes
the runs performed for the merged formulations of the methods  with 1, 4, and
16~RHS vectors on the Lomonosov and Lomonosov-2 supercomputers. In case of
preconditioned methods the identity operator is used as a preconditioner, which
is equal to single vector copy operation. While the parallelization inside the
node is implemented using MPI, the overhead due to communications for the chosen
test matrix is negligible compared the overall execution time, and results of
the comparison mostly validate correctness of the estimates for the
computational operations.

The obtained calculation results for all the methods considered in the paper
together with the predicted execution times are presented in
Tab.~\ref{tab:loop_exec_time_model_calc_Lom1} and
Tab.~\ref{tab:loop_exec_time_model_calc_Lom2} (results for the Lomonosov
supercomputer include only 1 and 4~RHS vector calculations due to compute node
memory limitations). The memory consumption to store the data for this test case
is measured by several gigabytes, which allows to suggest the use of memory
bandwidth of the random-access memory (RAM) as a data transfer bandwidth model
estimation. The presented results demonstrate good correspondence of the
theoretical and measured execution times for both compute systems. The
analytical model slightly underpredicts execution time for the Pipelined
BiCGStab method, providing the deviation by approximately 12\%. The difference
for all other test cases is within 5\%. The obtained predictivity is supposed to
be acceptable for the comparison purposes stated in the current paper.

\begin{table}[t!]
  \caption{Predicted and measured execution times of the single iteration loop
  for the BiCGStab methods with variable number of RHS vectors. Lomonosov
  supercomputer, times in seconds. \label{tab:loop_exec_time_model_calc_Lom1}}
\centering
 \begin{tabular}{ | c | c | c | c | c |}
 \hline
  \multirow{2}{*}{Method} & \multicolumn{2}{|c|}{$m = 1$}
  & \multicolumn{2}{|c|}{$m = 4$} \\
  \cline{2-5}
         & Model & Real & Model & Real \\
\hline
	BiCGStab & 0.106 & 0.107 & 0.3 & 0.29 \\
\hline
 	IBiCGStab & 0.109 & 0.113 & 0.32 & 0.32 \\
\hline
 	PipeBiCGStab & 0.12 & 0.127 & 0.36 & 0.38 \\
\hline
	PBiCGStab & 0.115 & 0.116 & 0.34 & 0.33 \\
\hline
 	RBiCGStab & 0.119 & 0.117 & 0.35 & 0.34 \\
\hline
 	PPipeBiCGStab & 0.143 & 0.156 & 0.45 & 0.49 \\
\hline
\end{tabular}
\end{table}

\begin{table}[h!]
  \caption{Predicted and measured execution times of the single iteration loop
  for the BiCGStab methods with variable number of RHS vectors. Lomonosov-2
  supercomputer, times in seconds. \label{tab:loop_exec_time_model_calc_Lom2}}
\centering
 \begin{tabular}{ | c | c | c | c | c | c | c |}
 \hline
  \multirow{2}{*}{Method} & \multicolumn{2}{|c|}{$m = 1$}
  & \multicolumn{2}{|c|}{$m = 4$} & \multicolumn{2}{|c|}{$m = 16$} \\
  \cline{2-7}
         & Model & Real & Model & Real & Model & Real \\
\hline
	BiCGStab & 0.057 & 0.058 & 0.16 & 0.16 & 0.59 & 0.58 \\
\hline
 	IBiCGStab & 0.059 & 0.061 & 0.17 & 0.18 & 0.62 & 0.64 \\
\hline
 	PipeBiCGStab & 0.065 & 0.069 & 0.19 & 0.21 & 0.71 & 0.76 \\
\hline
	PBiCGStab & 0.062 & 0.064 & 0.18 & 0.19 & 0.66 & 0.69 \\
\hline
 	RBiCGStab & 0.064 & 0.065 & 0.19 & 0.19 & 0.7 & 0.69 \\
\hline
 	PPipeBiCGStab & 0.077 & 0.085 & 0.24 & 0.27 & 0.9 & 1.01 \\
\hline
\end{tabular}
\end{table}

\subsubsection{Multi-node validation}

Following~\cite{Cools2017, Zhu2014} the multi-node validation is performed for
the test matrix, obtained as a result of discretization of 2D problem. The
matrix of 1~mln. unknowns, corresponding to the grid of $1000^2$ cells and
5-point stencil is used in the tests. The native ordering allows to obtain
almost constant communication overhead for the SpMV operation in the range up to
$10^3$ compute processes, thus allowing to evaluate the impact of global
communications on the overall methods execution time. The suboptimal data
decomposition for the general matrices can be obtained using the graph
partitioning algorithms (e.g.,~\cite{Chevalier2008, Karypis2003}). This,
however, lies beyond the scope of this paper.

\begin{figure}[tb]
\centering
\includegraphics[width=9cm]{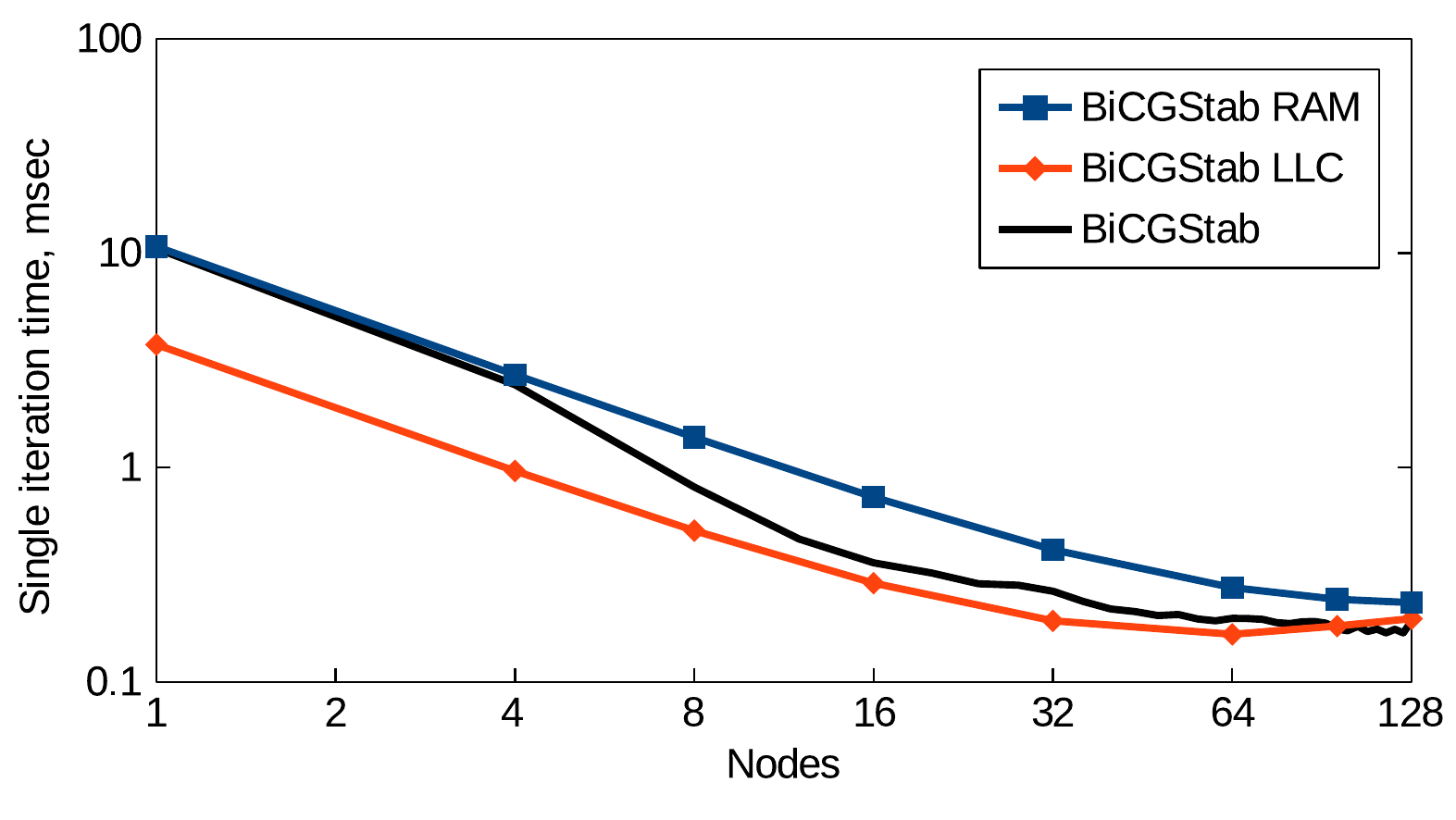}
\vspace{-0.5cm}
\caption{Analytical predictions and numerical simulation results for the
classical BiCGStab method.}
\label{fig:model_validation}
\end{figure}

Expressions~\eqref{est:bicgstab}-\eqref{est:ppipebicgstab} allow to predict the
expected execution times for the single iteration of the numerical methods
depending on the number of compute nodes and RHS vectors. The only ambiguous
parameter in the estimates is the memory bandwidth. Analytical expressions allow
to cover two limiting cases expecting the data is loaded from the RAM or the
data is loaded to the registers from the last level cache. The corresponding
predictions for the BiCGStab method with single RHS vector are presented in
Figure~\ref{fig:model_validation}. Starting with about 3.7~times difference for
the single node runs, these curves become closer with increasing the number of
compute nodes involved in the calculations. At the scale of $10^2$ nodes
(800~processes) the scalability reaches saturation, and predicted times become
almost identical for both configurations, which indicates domination of
communications over computations.

The measured execution times for the BiCGStab method fit the range between the
two analytical curves. Computations with 1--4 nodes provide results
corresponding to the RAM bandwidth predictions. The further increase in the number of compute
nodes allows to decrease the memory consumption per node and improve the cache
reuse. Starting from several tens of nodes the data fits the cache, which leads
to a reduction in the execution time and shift of the calculation results
towards the LLC~predictions.

\begin{figure}[tb] 
\centering 
\includegraphics[width=6cm]{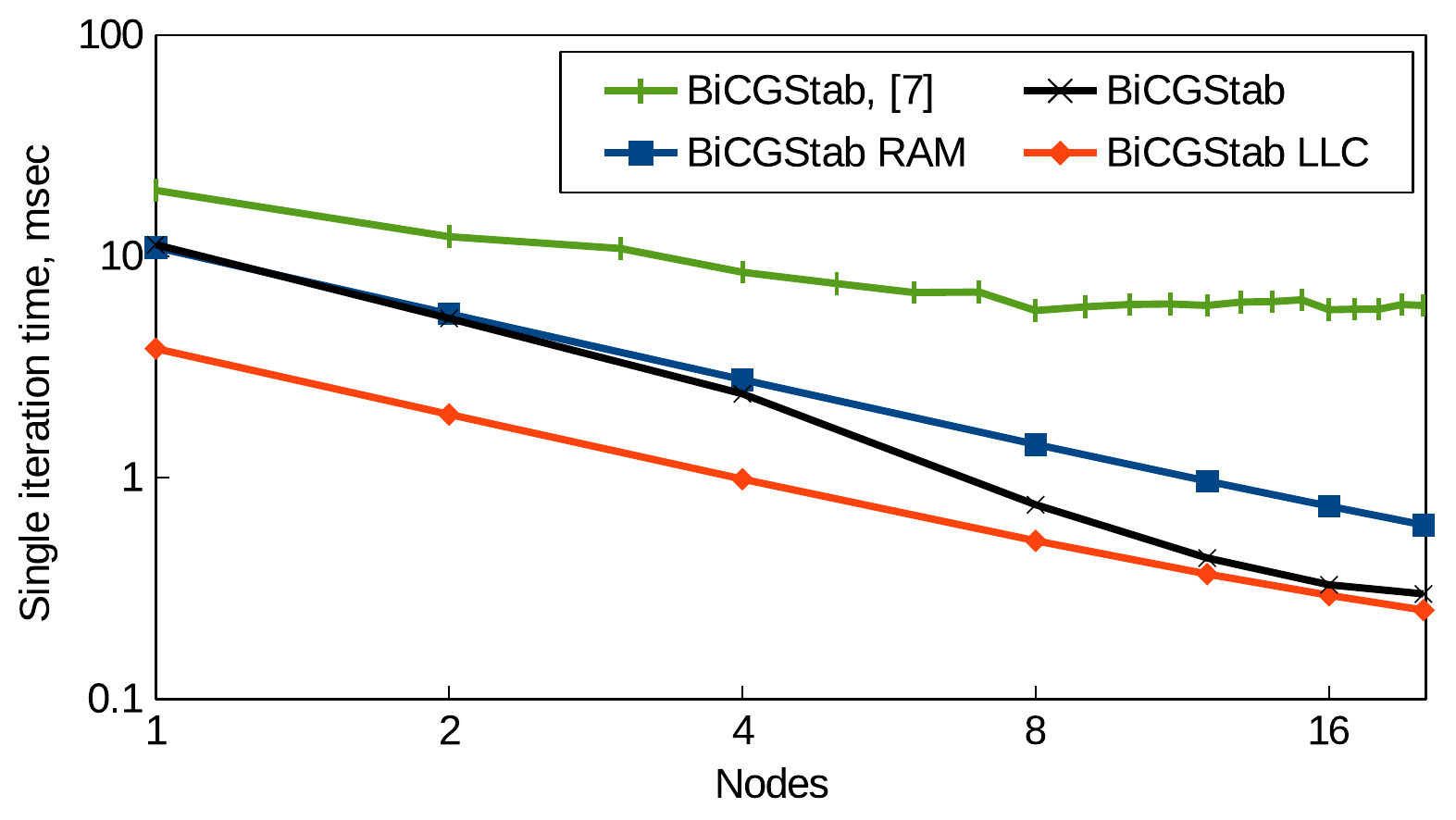}
\includegraphics[width=6cm]{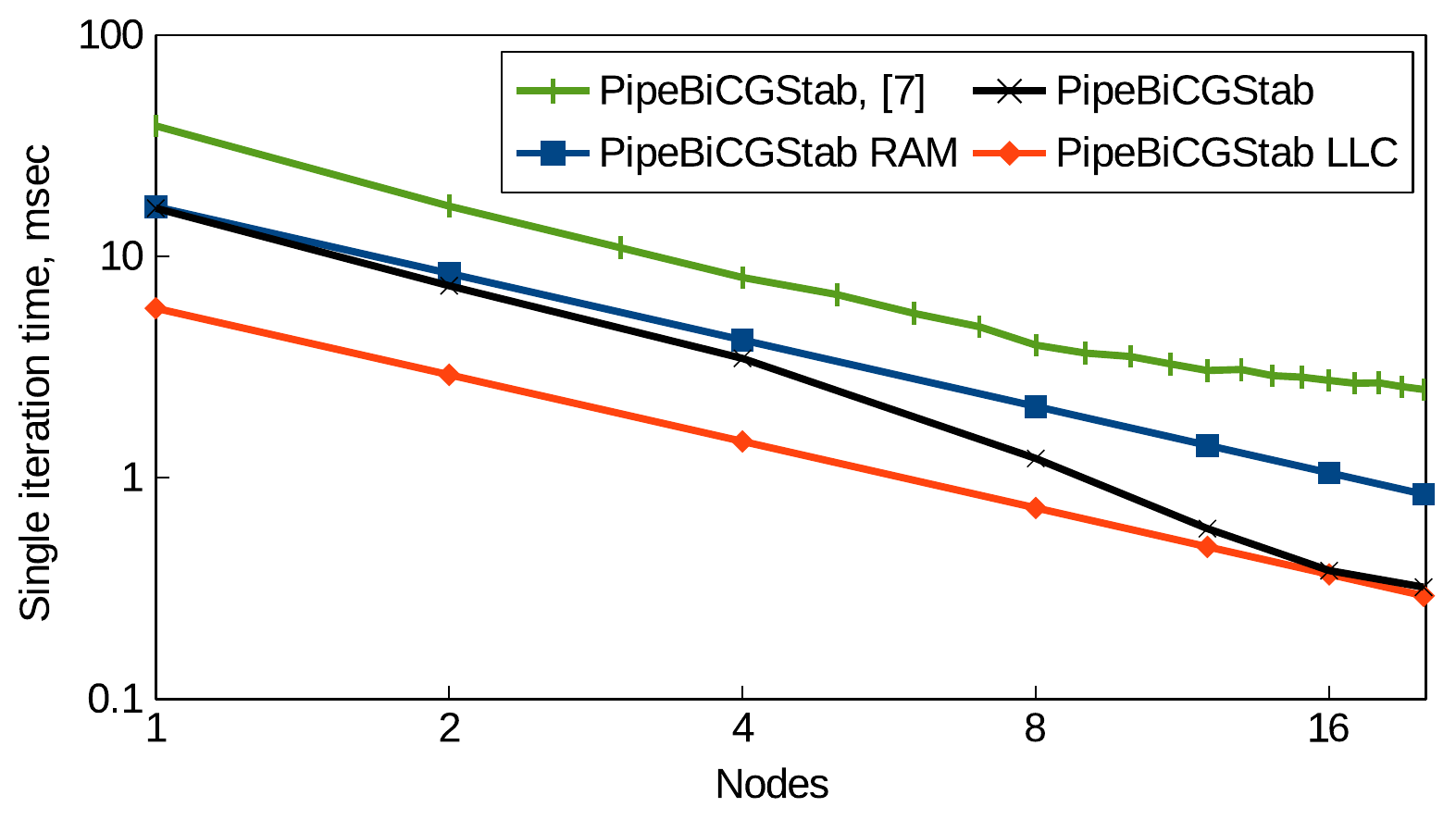}
\vspace{-1.0cm}
\caption{Comparison of execution times for the basic formulations of
classical BiCGStab and Pipelined BiCGStab methods with results presented
in~\cite{Cools2017}.}
\label{fig:cools_timings_comparison}
\end{figure}

The paper~\cite{Cools2017} contains execution time results for 1~mln. 5-diagonal
test matrix, performed on the hardware platform very similar to the one used for
the validation in the current paper (InfiniBand QDR interconnect, compute nodes
with 2 x 6-core Intel Xeon X5660 processors) with only minor difference in the
number of processor cores per node. To compare the model predictions and
simulation results together with results of other authors the execution times
for the single iteration of basic formulations of classical BiCGStab and
Pipelined BiCGStab methods are summarized in
Figure~\ref{fig:cools_timings_comparison}. The present simulation results
demonstrate good agreement with analytical model predictions for both iterative
methods. The obtained results significantly outperform the ones presented
in~\cite{Cools2017}. The newly developed implementation of the methods provides
1.8 and 2.4~times lower execution times for the single-node runs of classical
BiCGStab and Pipelined BiCGStab methods correspondingly.

\begin{figure}[tb] 
\centering
\includegraphics[width=6cm]{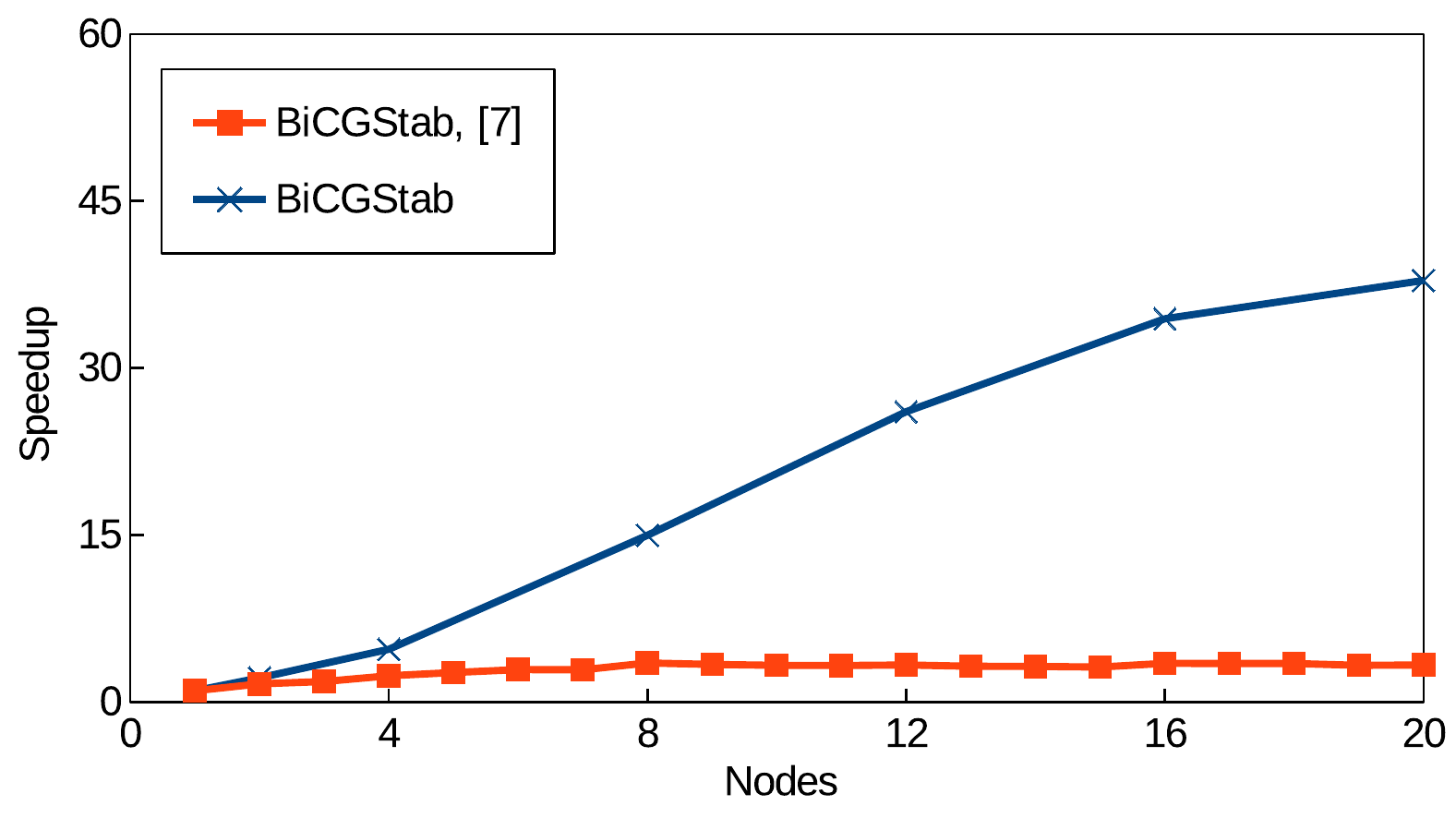}
\includegraphics[width=6cm]{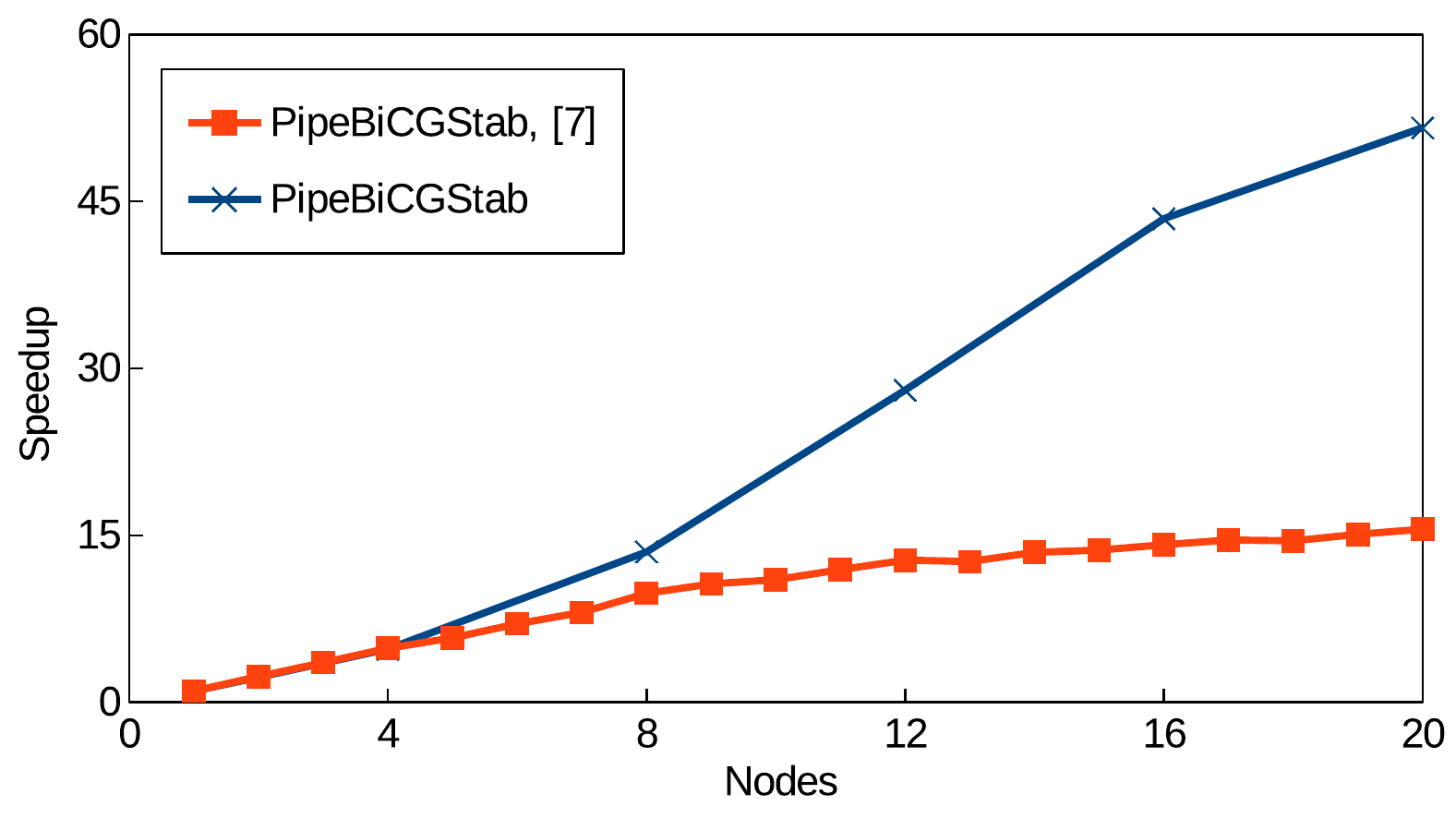}
\vspace{-1.0cm}
\caption{Scalability of basic formulations of classical BiCGStab and Pipelined
BiCGStab methods compared with results presented in~\cite{Cools2017}.}
\label{fig:cools_scalability_comparison}
\end{figure}

The poor scalability of the methods for the considered test matrix is
demonstrated in~\cite{Cools2017}. The speedup, defined as a ratio of the
execution time on the single node to the one with $p$ nodes,
\begin{gather}
S(p) = \frac{T(1)}{T(p)}.
\end{gather}
does not exceed~3.5 for the classical BiCGStab, and for the Pipelined BiCGStab
it reaches~15.5. The present results, to the contrary, demonstrate the
superlinear speedup (Figure~\ref{fig:cools_scalability_comparison}). This effect
is a result of cache usage improvement, which allowed more than 35-fold
decrease of the execution time on 20~compute nodes. In total, the current
implementation of the iterative methods outperforms results of~\cite{Cools2017}
for the classical BiCGStab by a factor of~20, and for the Pipelined BiCGStab by
a factor of~7.8. While the paper~\cite{Cools2017} shows advantage of the
Pipelined BiCGStab method starting from 4~compute nodes, the present results do
not reproduce this behaviour. The classical method outperforms the pipelined
variant in the whole range of the compute nodes considered in this test.

%%%%%%%%%%%%%%%%%%%%%%%%%%%%%%%%%%%%%%%

\section{Comparison of modified BiCGStab methods}

The above presented results of the analytical model validation demonstrate good
correspondence of predictions and calculations. This allows to use the proposed
model to investigate the influence of various parameters on the execution times
of the methods and indicate the range of applicability for each of the iterative
methods considered in the paper.

\subsection{Model predictions}

The range of applicability for the modified BiCGStab methods can be estimated by
analysing expressions~\eqref{est:bicgstab}-\eqref{est:ppipebicgstab}. The
modified methods allow to decrease the number of global reductions by the price
of several extra vector operations and hide the communications by computations.
The typical scale of the problem per each compute node when the modified method
can outperform the classical one can be estimated as:
\begin{gather}
T_G \sim \frac{8 N m}{b p}.
\end{gather}
Using the reference values for global reduction of 20 $\mu$s and LLC memory
bandwidth of $10^{11}$ bytes/s, one can obtain the estimate $N m/p \sim 2.5
\cdot 10^4$ elements per node, and this value further decreases with increasing
the number of RHS vectors.

The further analysis of the range of optimality for the numerical methods is
performed with help of the following parameter:
\begin{gather}
R^i(p) = \frac{\min_{j \in K} \left( T^j(p) \right)}{T^i(p)}, \, i \in K,
\end{gather}
which can be interpreted as a relative performance of the specific method among
the set of the methods considered. Here, $K$ is the set of the iterative
methods, i.e. \textit{K = \{BiCGStab, IBiCGStab, PipeBiCGStab\} } or \textit{K =
\{PBiCGStab, RBiCGStab, PPipeBiCGStab\}}. This parameter varies in the range
$R^i \in (0,1]$ and shows how close is the specific method to the optimal
numerical method providing the minimal execution time. The proposed
estimates~\eqref{est:bicgstab}-\eqref{est:ppipebicgstab} allow to plot the
corresponding distributions depending on various model parameters as a functions
of the number of compute nodes. These plots clearly show the range of optimality for
each of the considered methods.

\begin{figure}[tb] 
\centering
\includegraphics[width=12cm]{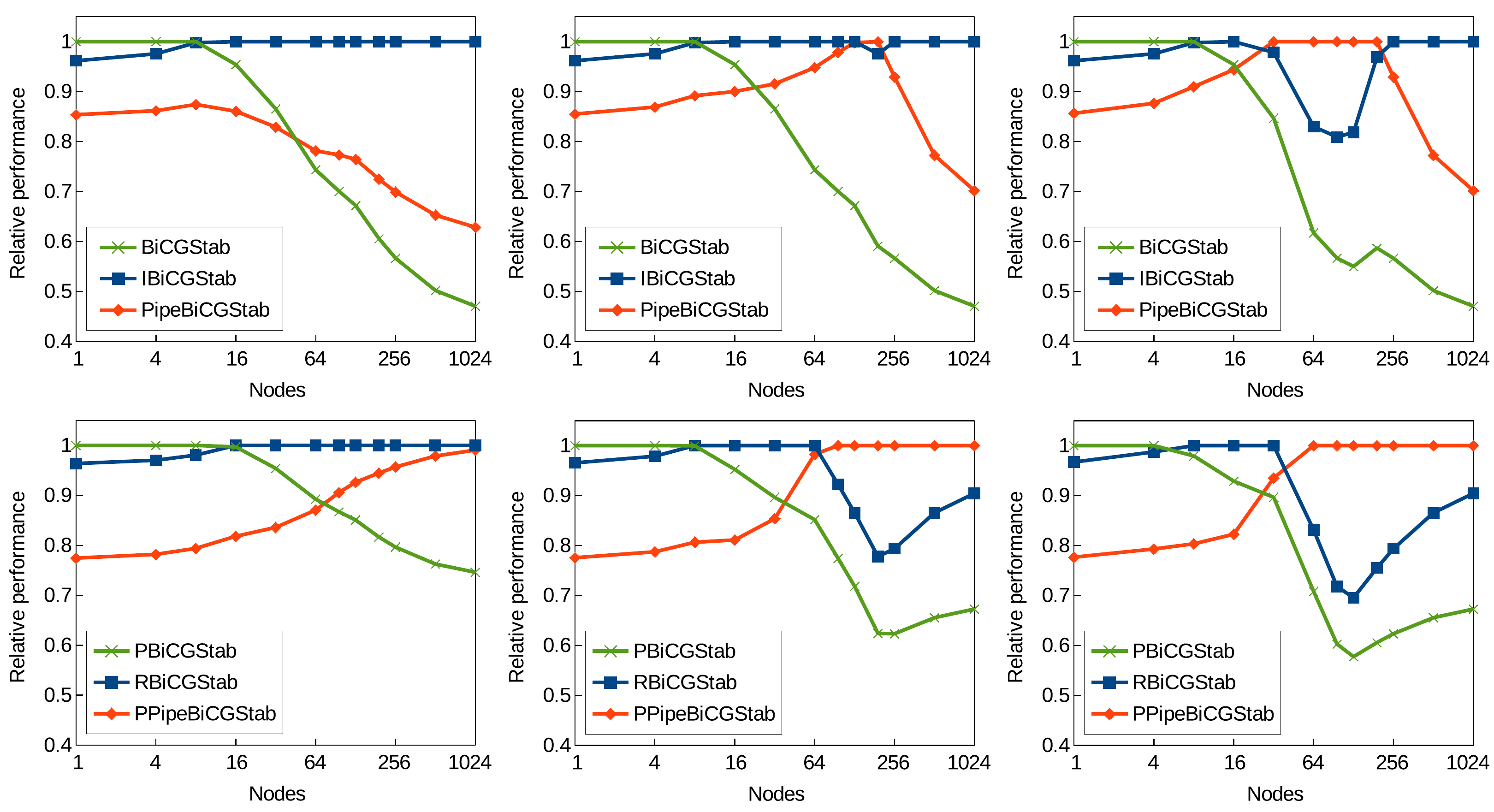}
\vspace{-0.5cm}
\caption{Relative performance of iterative methods when solving SLAE with single
RHS vector depending on degree of overlap of communications and calculations. From
left to right: no overlap, $\gamma = 1$; partial overlap, $\gamma = 0.5$; full
overlap, $\gamma = 0$.}
\label{fig:model_prediction_m1_gamma}
\end{figure}

The plots, presented in Figure~\ref{fig:model_prediction_m1_gamma}, allow to
compare the efficiency of the numerical methods when solving system of linear
algebraic equations with single RHS vector, and indicate the influence of
overlapping of global communications with computations. The lowest execution
times at the scale of 1--10~compute nodes among the unpreconditioned methods are
observed for the classical BiCGStab method. The Improved BiCGStab has about 5\%
penalty and the Pipelined BiCGStab demonstrates about 15\% relative performance
degradation. The Improved BiCGStab outperforms two other methods at the scale of
10--30~compute nodes due to only single global reduction performed at each
iteration compared to three global reductions needed by the classical BiCGStab
and two global reductions of the pipelined method. The choice of the optimal method
at the higher scales of compute nodes varies depending on the overlap parameter
$\gamma$. In case of no overlap ($\gamma = 1$) the Improved BiCGStab remains the
optimal method in the range of compute nodes considered. The partial overlap
($\gamma=0.5$) improves the characteristics of the Pipelined BiCGStab method,
and at the scale of 100 nodes it provides the execution times compatible with
the Improved BiCGStab. Finally, the full overlap ($\gamma = 0$) allows to
unleash the potential of the Pipelined BiCGStab method. This method provides the
highest performance in the range of 30--250 compute nodes.

In case of using the lightweight preconditioner (the computational complexity of
the preconditioner, $\alpha$, equal to two vector reads/writes) the classical
BiCGStab is also the optimal method at the scales of 1--10 nodes. The Reordered
BiCGStab outperforms the classical BiCGStab method at the scales of
10--50~compute nodes due to reduced number of global reductions. The
preconditioned Pipelined BiCGStab has the same number of global reductions, but
significantly higher number of vector operations, thus demonstrating lower
performance compared the Reordered BiCGStab. In case of at least partial overlap
of global communications by computations the pipelined method can outperform two
other methods starting from 50--100~compute nodes. Its performance, however,
strongly depends on the overlap parameter,~$\gamma$.

\begin{figure}[tb] 
\centering
\includegraphics[width=12cm]{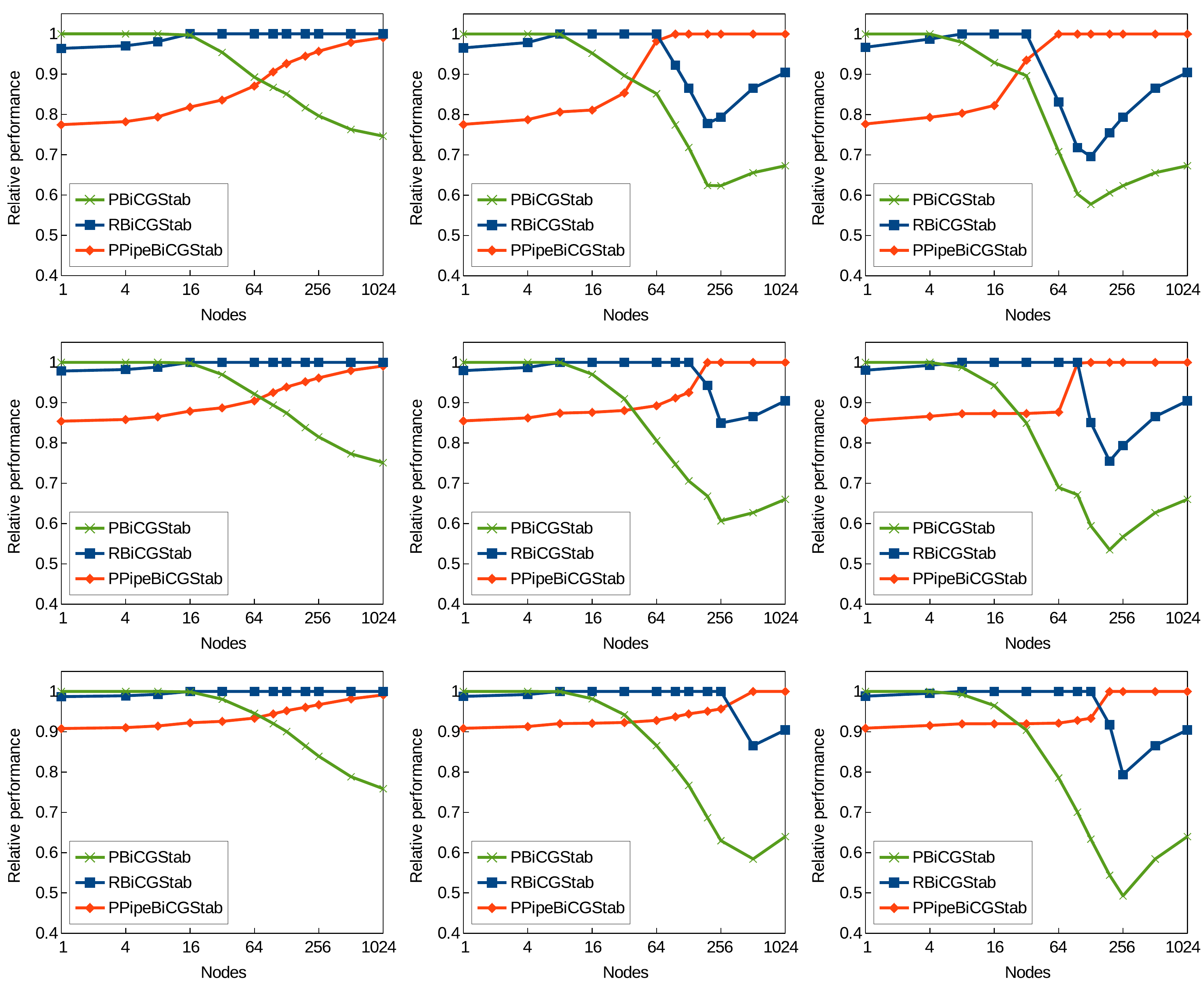}
\vspace{-0.5cm}
\caption{Influence of preconditioner complexity on the performance of iterative
methods. From left to right: no overlap, $\gamma = 1$; partial overlap, $\gamma
= 0.5$; full overlap, $\gamma = 0$. From top to bottom: no preconditioner,
$\alpha = 2$; $\alpha = 20$; $\alpha = 50$.}
\label{fig:model_precond_m1_gamma}
\end{figure}

The preconditioner complexity does not have principal influence on the methods
relative performance (Figure~\ref{fig:model_precond_m1_gamma}). The increase of
the preconditioner complexity shifts towards the higher scales the transition
point from the Reordered BiCGStab to the Pipelined BiCGStab and reduces the
difference in the relative performance of the methods. The latter one is caused
by a reduction of the role of extra vector operations in the overall
computational costs.

\begin{figure}[tb!] 
\centering
\includegraphics[width=12cm]{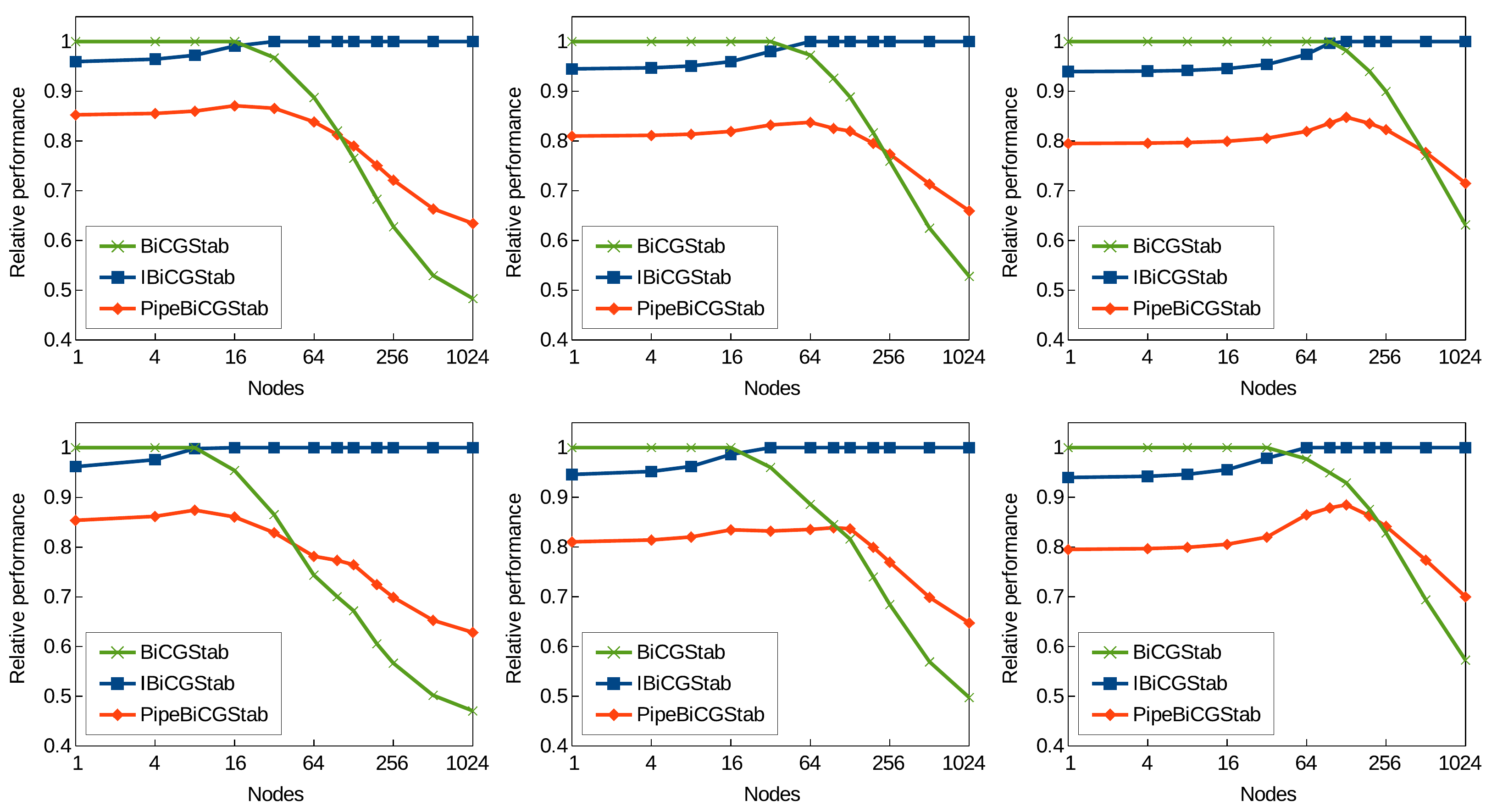}
\vspace{-0.5cm}
\caption{Relative performance of the iterative methods for SLAEs with multiple
RHS vectors, unpreconditioned methods. First row -- RAM memory bandwidth
estimates; second row -- LLC memory bandwidth estimates. From left to right: $m
= 1$; $m = 4$; $m = 16$.}
\label{fig:model_mrhs_unprec}
\end{figure}

\begin{figure}[tb!] 
\centering
\includegraphics[width=12cm]{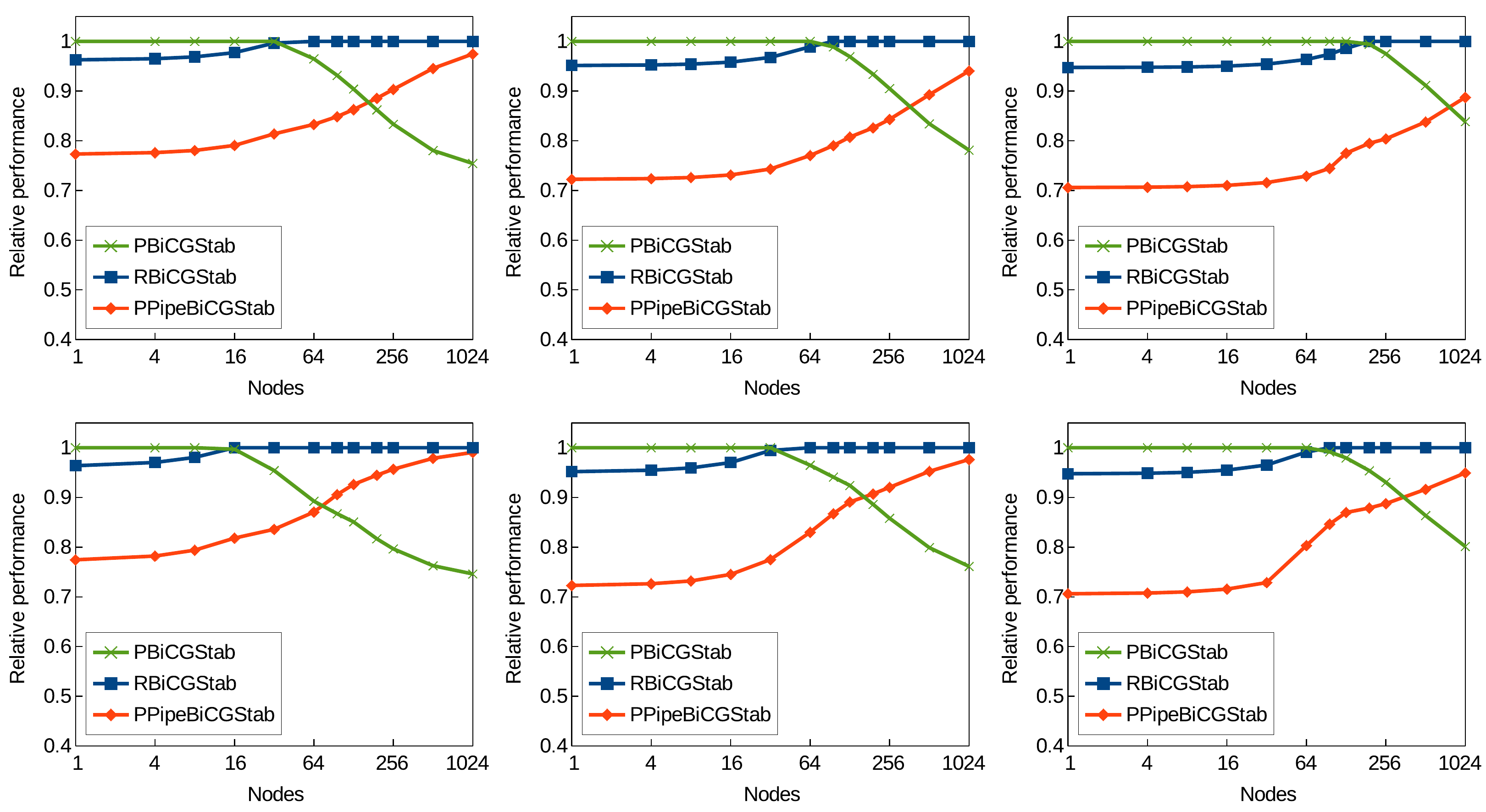}
\vspace{-0.5cm}
\caption{Relative performance of the iterative methods for SLAEs with multiple
RHS vectors, preconditioned methods. First row -- RAM memory bandwidth
estimates; second row -- LLC memory bandwidth estimates. From left to right: $m
= 1$; $m = 4$; $m = 16$.}
\label{fig:model_mrhs_prec}
\end{figure}

The increase in the number of RHS vectors increases the memory consumption and
a fraction of vector operations in the overall execution time. The first aspect
leads to the increase in the number of compute nodes when the data starts to fit
the cache. While for the single RHS vector the cumulative LLC memory of about
15~compute nodes of Lomonosov supercomputer is enough to store the data, for 4
and 16~RHS vectors this value increases to 30~and 100~nodes correspondingly.
This indicates that the expected transition points for the simulations with
4~RHS vectors would be in between of analytical model predictions with RAM and
LLC memory bandwidth, and results for 16~RHS vectors would mostly correspond to
the RAM bandwidth values.

The analytical model predictions for 1, 4, and 16~RHS vectors with RAM and LLC
memory bandwidths are summarized in
Figures~\ref{fig:model_mrhs_unprec}-\ref{fig:model_mrhs_prec}. The plots
indicate that the increase in the amount of RHS vectors for the unpreconditioned
methods (Figure~\ref{fig:model_mrhs_unprec}) leads to a systematic shift of the
transition point from the classical BiCGStab to the Improved BiCGStab towards
the higher scales in compute nodes. The Pipelined BiCGStab demonstrates at least
20\% higher execution times compared to the other methods in the whole range of
compute nodes considered.

The similar situation is observed for the preconditioned methods
(Figure~\ref{fig:model_mrhs_prec}): the preconditioned Pipelined BiCGStab
demonstrates about 30\% higher execution times and this value decreases only at
the scale of thousands of compute nodes. The Reordered BiCGStab method can be
beneficial compared the preconditioned BiCGStab method starting from the scales
of 32--128~compute nodes, and this transition point also shifts to the right
along the compute nodes axis with increasing the number of RHS vectors.

The proposed analytical model allowed to compare execution times of the BiCGStab
methods and identify the influence of various parameters. The estimates presented above
cover the range of compute nodes from 1 to 1024. Accounting the size of the test
problem and the specific form of the local communications time function, which
is generally valid up to 1000 compute processes (about 120~compute nodes of
the current test platform), it should be noted that practical range of these
estimates is limited by 128~nodes. To be able to predict correctly the execution
times at the higher scales, the more complicated local communication time
function must be applied.

%%%%%%%%%%%%%%%%%%%%%%%%%%%%%%%%%%%%%%%

\subsection{Validation of model predictions with simulation results}

The current section focuses on validation of results obtained using the proposed
analytical execution time model. The validation starts with investigation of the
single iteration execution times when solving SLAEs with 1, 4, and 16~RHS
vectors for the methods considered in the paper. The corresponding experiments
are performed for the merged formulations of the methods in the range of
1--128~compute nodes and utilize all available 8~cores per each node. This test
series do not assume any special manipulations with progression of non-blocking global
communications. Accounting results of the corresponding benchmark, presented
in~\ref{app:Iallreduce_async}, no overlap of communications is expected.

\begin{figure}[tb] 
\centering
\includegraphics[width=12cm]{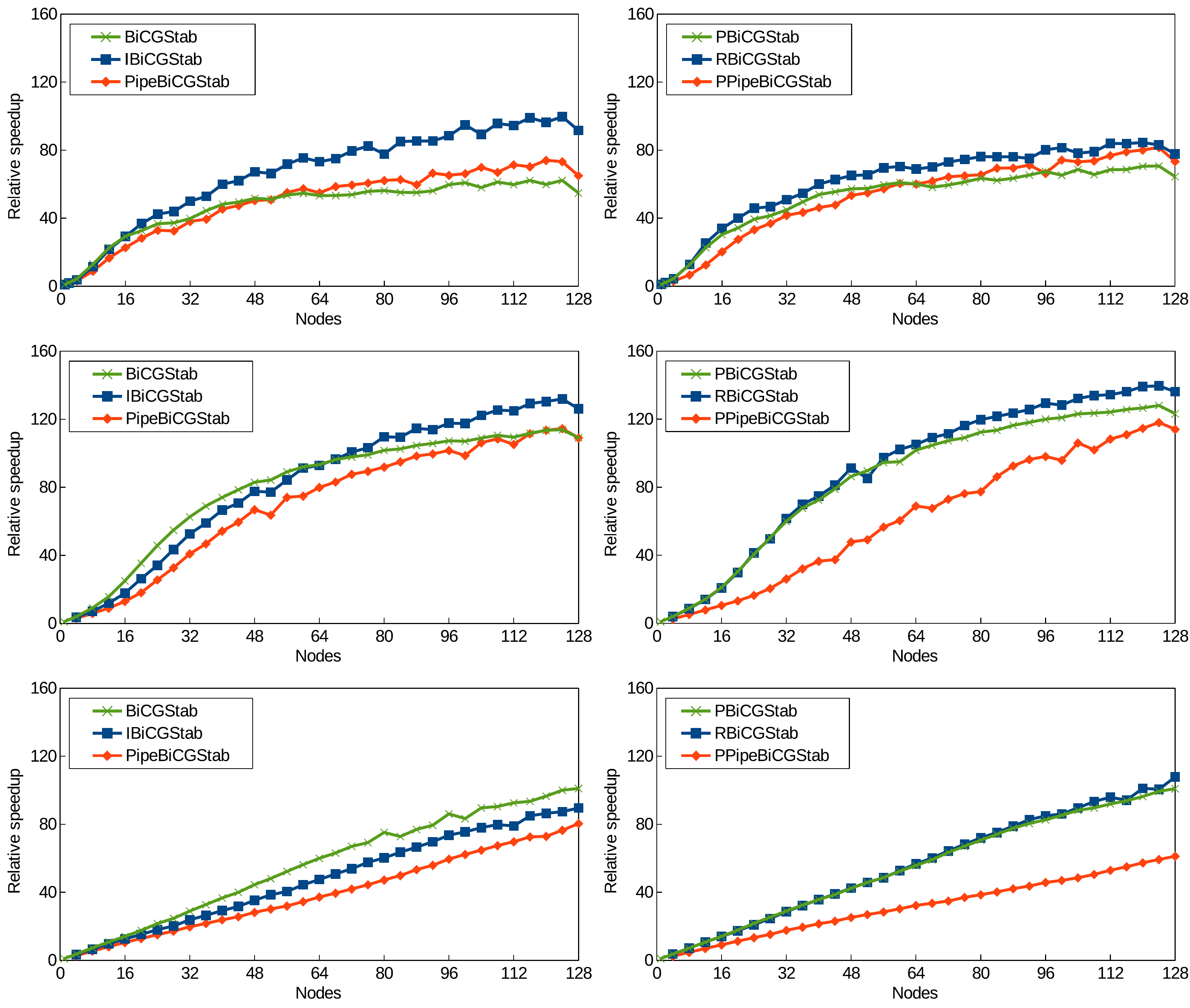}
\vspace{-0.5cm}
\caption{Relative scalability of the methods when solving SLAEs with
multiple RHS vectors. From top to bottom: $m=1$; $m = 4$; $m = 16$.}
\label{fig:sim_mrhs}
\end{figure}

The measured simulation results in the form of relative speedup are presented in
Figure~\ref{fig:sim_mrhs}. This parameter is defined as a ratio of the execution
time for the classical BiCGStab method (or the preconditioned one) when
performing simulation on the single compute node to the execution time of the
specific method with $p$ compute nodes:
\begin{gather}
P^i(p) = \frac{T^{(P)BiCGStab}(1)}{T^i(p)}.
\end{gather}
The presented plots demonstrate good correspondence with the analytical model
predictions for the case with no overlap of global communications ($\gamma =
1$). The analytical model correctly predicts the scales of compute nodes
corresponding to the transition points. Simulations with single RHS vector
confirm predictions in that the optimal methods for the scales of 1--10~compute
nodes is the classical BiCGStab, and for the higher scales -- Improved BiCGStab
or Reordered BiCGStab among the unpreconditioned and preconditioned methods
correspondingly. For the unpreconditioned methods the Improved BiCGStab becomes
the optimal method for the simulations with 4~RHS vectors starting from about
64~nodes while for the lower scales the classical BiCGStab demonstrates the best
performance. The further increase of the amount of RHS vectors increases the
role of extra vector operations and moves the transition point between the
classical BiCGStab and the Improved BiCGStab beyond the scale of 128~nodes.
While the Pipelined BiCGStab has high potential for large scale simulations, the
overlap of global communications by computations is a prerequisite to unleash
all the advantages of this method. The absence of overlap of global
communications, however, does not allow to demonstrate advantages of the
Pipelined BiCGStab method in practice. 

The preconditioned BiCGStab and Reordered BiCGStab methods demonstrate about the
same relative scalability results, with minor advantage of classical method at
the lower scales and of Reordered BiCGStab at the higher scales. The performance
variation, however, is only within several percent. The preconditioned
Pipelined BiCGStab provides significantly higher execution times, which is a
result of higher number of extra vector operations occurred as a result of
algorithm reordering.

\begin{figure}[tb]
\centering
\includegraphics[width=12cm]{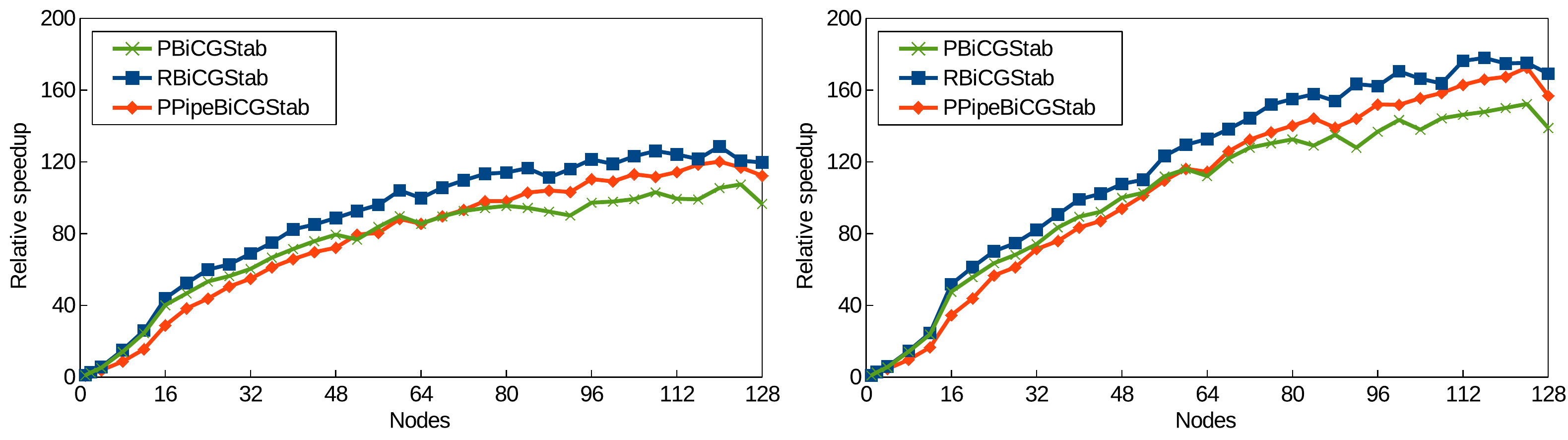}
\vspace{-0.5cm}
\caption{Relative scalability of the methods with various preconditioner
complexity when solving the SLAE with single RHS vector. Left -- $\alpha=20$;
right -- $\alpha = 50$.}
\label{fig:sim_m1_precond}
\end{figure}

The second validation case investigates the influence of the preconditioner
complexity on the relative speedup of the methods. The corresponding runs are
performed for the group of preconditioned methods and two preconditioners with
complexity equal to 20 and 50~vector read/write operations. The obtained
scalability results demonstrate similar behaviour for both configurations with
only the difference in the maximal values (Figure~\ref{fig:sim_m1_precond}). The
increase in the preconditioner complexity leads to an increase in the peak
scalability values by a factor of~1.5. The methods optimality intervals
correspond to the ones obtained as a result of analytical model predictions.

Presented above simulation results and results in~\ref{app:Iallreduce_async}
demonstrate that progression of global non-blocking operations is not
automatically performed during the calculations on the compute system used in
the parallel runs. Lomonosov supercomputer has no hardware support for the
non-blocking collective operations, and the use of software progression,
implemented in the Intel MPI library with specialized progression threads, leads
to significant performance degradation, even though it provides the expected
communications overlap. This fact makes corresponding functionality inapplicable
for the real simulations. As an alternative, the manual progression is
investigated in the following series of numerical experiments. The scalability
for the merged formulations of Reordered BiCGStab and preconditioned Pipelined
BiCGStab methods with manual progression is investigated in the range of
1--128~compute nodes. The preconditioner is emulated in the form of a loop over
a set of vector operations, and the \textit{MPI\_Test} calls are performed after
each vector operation. The preconditioners with complexity of 20 and 50 vector
read/write operations are considered, i.e. 10 and 25 \textit{MPI\_Test} calls
can be performed to progress global reduction. For the pipelined method
additional \textit{MPI\_Test} calls are also performed during the SpMV
operation.

\begin{figure}[tb]
\centering
\includegraphics[width=12cm]{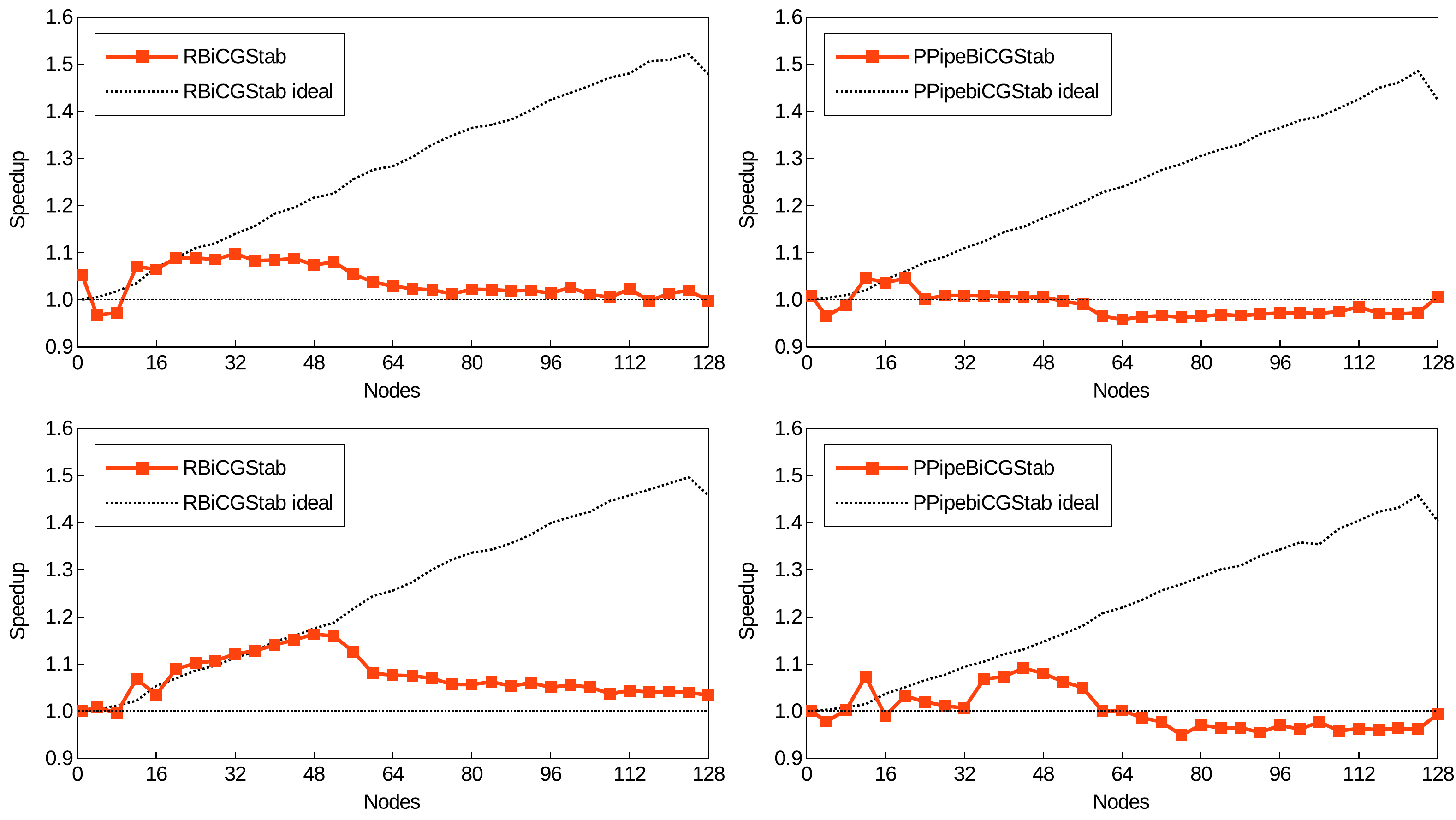}
\vspace{-0.5cm}
\caption{Simulation speedup due to manual progression of global non-blocking
collective operations for preconditioners of different complexity. Top --
$\alpha = 20$; bottom -- $\alpha = 50$.}
\label{fig:sim_progression}
\end{figure}

The obtained calculation results as a reference of execution times for the basic
implementation to the one with manual progression are presented in
Figure~\ref{fig:sim_progression}. The values, greater than~1, indicate the
speedup of the simulations due to manual progression. These plots also contain
the approximate values, showing the expected performance for the case of ideal
overlap of global communications with computations (here the expected execution
time is defined as the measured execution time of basic implementation minus
expected global communication time~\eqref{eq:global_comm}). The
Figure~\ref{fig:sim_progression} demonstrates that the manual progression allows
to obtain the simulation speedup and achieve the communications overlap in real
simulations, but its efficiency depends on lots of factors. The case with lower
preconditioner complexity (and lower number of \textit{MPI\_Test} calls)
demonstrates observable speedup for the Reordered BiCGStab in the range of
12--64~compute nodes. The speedup for 12--24~nodes corresponds to the ideal
overlap, but further increase in the compute nodes shows the degradation of
overlap efficiency. 
The effect of decreasing the overlap efficiency may be caused by the amount of
computations in the preconditioner becoming too small, and the effect seems to
be increasing with growing the number of compute nodes.
The observed speedup due to manual
progression does not exceed 10\%, and tends to zero when using more than
64~nodes. The preconditioned Pipelined BiCGStab method demonstrates lower
efficiency of overlapping. Manual progression allows to speedup the calculations
for the range of 12--20~compute nodes and the speedup of about 5\% is achieved.
The use of 24--56~nodes provides the same execution times as the basic
implementation, and for the higher scales the slowdown by 3-5\% is observed.

The increase in the preconditioner complexity (and in the number of
\textit{MPI\_Test} calls) improves the communications overlap and allows to
obtain the speedup for the Reordered BiCGStab of about~20\%. Range of the nodes
with ideal communications overlap is increased up to~48, and the 5\% speedup is
observed for the higher scales. The similar tendency is reproduced for the
Pipelined BiCGStab method, but the peak speedup is limited by 10\%, and the use
of progression for more than 48~nodes still leads to a slowdown by 2-5\%.

%%%%%%%%%%%%%%%%%%%%%%%%%%%%%%%%%%%%%%%%%%%%%%%%%%%%%%%%%%%%

\section{Conclusions}

The efficiency of the classical BiCGStab iterative method and several modified
formulations for solving systems of linear algebraic equations is investigated.
The importance of reducing the costs of vector operations to achieve maximal
performance is demonstrated. The vector merging technique is applied for the
methods considered, and corresponding merged formulations of the methods
allowing to minimize data transfers with the memory are proposed. While the
merging of vector operations is applicable for all the methods, the effect is
much more substantial for the modified formulations. The suggested optimization
allows to reduce the extra vector operation costs and decrease the gap between
the classical method and modified formulations. The observed speedup due to
merging of vector operations is up to~30\%.

The analytical execution time model is proposed to perform the detailed
comparison of the methods. The model is based on the volume of data transfers
with the memory. The accuracy of this model is demonstrated by the validation
with the simulation results. The analytical model allows to predict the
corresponding execution times with about 10\%~error.

The proposed analytical model is used for the detailed comparison of six
iterative methods. The influence of the problem size, preconditioner complexity,
number of RHS vectors, and the effect of global communications overlap on the
choice of the optimal method are highlighted. It is shown that the modified
methods can outperform the classical BiCGStab method when the global reduction
time becomes comparable with the vector update time. The increase in the
preconditioner complexity and the number of RHS vectors extends the range of
optimality for the classical BiCGStab methods towards the higher scales in
compute nodes. The corresponding theoretical predictions are validated by the
numerical simulations, and the simulation results completely reproduce the
predicted behaviour.

The progression of non-blocking global communications is among the key features
required by the modified BiCGStab methods. In practice, however, the efficient
asynchronous progression for the global reductions with short messages
(8--8192~bytes) is a challenging issue. The manual progression is investigated
and implemented in the corresponding test code. The performed experiments have
shown that the manual progression can generally produce the overlap of
communications by computations and speedup the calculations, but is problematic
to be implemented in the computational codes.

%%%%%%%%%%%%%%%%%%%%%%%%%%%%%%%%%%%%

\appendix

\section{Compute platforms}
\label{section:spec}

The results of numerical experiments discussed in the paper were performed on
compute systems Lomonosov and Lomonosov-2 installed in the Research Computing
Center of Lomonosov Moscow State University. The main parameters of the
supercomputers are summarized in tab.~\ref{tab:spec}. The table also contains
the real memory bandwidth values, obtained with help of STREAM
benchmark~\cite{McCalpin1995}. It should be noted the last level cache (LLC)
bandwidth values demonstrate observable deviations from run to run, thus the
presented results should be treated as some reference values only.

\begin{table}[t!]
  \caption{The main characteristics of Lomonosov and Lomonosov-2
  supercomputers.\label{tab:spec}}
\centering
\begin{tabular}{ | c | c | c | c |}
\hline
Supercomputer & Lomonosov & Lomonosov-2 \\
\hline
	Processor model & Intel Xeon X5570 & Intel Xeon E5-2697v3 \\
\hline
	Processors & 2 & 1 \\
\hline
	Cores & 4 & 14 \\
\hline
	Instruction set & SSE4.2 & AVX2 \\
\hline
	LLC size, MB & 8 & 35 \\
\hline
	RAM Bandwidth, GB/s & 16 & 60 \\
\hline
	LLC Bandwidth, GB/s & 46 & 288 \\
\hline
	Interconnect & Infiniband QDR & Infiniband FDR \\
\hline
\end{tabular}
\end{table}

%%%%%%%%%%%%%%%%%%%%%%%%%%%%%%%%%%%%%%%%%%%%%%%%%%%%%%%%%%%

\section {Merged formulations of BiCGStab methods}
\label{app:merged_methods}

The current section presents merged formulations of the modified BiCGStab
methods discussed in the paper. The operations specified in the single line are
performed in the merged form, i.e. with single loop over the vectors.

\begin{minipage}[t]{0.46\textwidth}
\begin{algorithm}[H]
\small
    \centering
    \caption{The BiCGStab method~\cite{Vorst1992}}\label{Alg:BiCGStab}
    \begin{algorithmic}[1] \STATE $x_{0}$ = initial guess; $r_{0} = b - A x_{0}$
	  \STATE $\rho_{0} = (r_{0},r_{0})$
	  \STATE $p_{0} = r_{0}$\\
	  \FOR {$j=0,1,\ldots$}
	  \STATE $v_j = A p_{j}$
	  \STATE $\delta_{j} = (v_j,r_{0})$
	  \STATE $\alpha_{j} = \frac{\rho_{j}}{\delta_j}$
	  \STATE $s_{j} = r_{j} - {\alpha_{j}} v_j$
	  \STATE $t_j = A s_{j}$
	  \STATE $\phi_{j} = (t_j,s_{j})$ \label{BiCGStab_dot_merge_1}
	  \STATE $\psi_{j} = {(t_j,t_j)}$ \label{BiCGStab_dot_merge_2}
	  \STATE $\omega_{j} = \frac{\phi_j} {\psi_j}$
	  \STATE $x_{j+1} = x_{j} + \alpha_{j} p_{j} + {\omega_{j}} s_{j}$
	  \STATE $r_{j+1} = s_{j} - {\omega_{j}} t_j$
	  \IF{$(r_{j+1}, r_{j+1}) < \varepsilon^2$} \label{BiCGStab_conv_check} 
	  \STATE \bf{break}
	  \ENDIF
	  \STATE $\rho_{j+1} = (r_{j+1},r_{0})$
	  \STATE $\beta_{j} = {\frac{\rho_{j+1}} {\rho_{j}}}
	  {\frac{\alpha_{j}} {\omega_{j}}}$
	  \STATE $p_{j+1} = r_{j+1} + {\beta_{j}} (p_{j}-{\omega_{j}} v_j)$
	  \ENDFOR
    \end{algorithmic}
\end{algorithm}
\end{minipage}
\hfill
\begin{minipage}[t]{0.46\textwidth}
\begin{algorithm}[H]
\small
    \centering
    \caption{The merged BiCGStab method}\label{Alg:MBiCGStab}
    \begin{algorithmic}[1]
	  \STATE $x_{0}$ = initial guess; $z = A x_{0}$
	  \STATE $r_{0} = b - z$, $\rho_{0} = (r_{0},r_{0})$, $p_{0} = r_{0}$
	  \FOR {$j=0,1,\ldots$}
	  \STATE $v_j = A p_{j}$
	  \STATE $\delta_{j} = (v_j,r_{0})$
	  \STATE $\alpha_{j} = \frac{\rho_{j}}{\delta_j}$
	  \STATE $s_{j} = r_{j} - {\alpha_{j}} v_j$
	  \STATE $t_j = A s_{j}$
	  \STATE $\phi_j = (t_j,s_{j})$, $\psi_j = (t_j,t_j)$, \\
	         $\theta_j = (s_j,s_j)$ \label{MBiCGStab_group1}
	  \STATE $\omega_{j} = \frac{\phi_j}{\psi_j}$
	  \IF{$(\theta_j - \omega_j \phi_j) < \varepsilon^2$} 
	    \STATE $x_{j+1} = x_{j} + \alpha_{j} p_{j} + {\omega_{j}} s_{j}$,\\
	           $r_{j+1} = s_{j} - {\omega_{j}} t_j$ \label{MBiCGStab_group2}
	    \STATE \bf{break}
	  \ENDIF 
	  \STATE $x_{j+1} = x_{j} + \alpha_{j} p_{j} + {\omega_{j}} s_{j}$, \\
	  		 $r_{j+1} = s_{j} - {\omega_{j}} t_j$, \\
	  		 $\rho_{j+1} = (r_{j+1},r_{0})$ \label{MBiCGStab_group3}
	  \STATE $\beta_{j} = {\frac{\rho_{j+1}} {\rho_{j}}}
	  {\frac{\alpha_{j}} {\omega_{j}}}$
	  \STATE $p_{j+1} = r_{j+1} + {\beta_{j}} (p_{j}-{\omega_{j}} v_j)$
	  \ENDFOR
    \end{algorithmic}
\end{algorithm}
\end{minipage}

%%%%%%%%%%%%%%%%%%%%%%%%%%%%%%%%%%%%%%%%%%%%%%%%%

\begin{algorithm}[H]
\small
  \centering
  \caption{The merged Improved BiCGStab method $^4$
%  \footnotemark
}
  \label{Alg:IBiCGStab}
  \begin{algorithmic}[1]
    \STATE $x_{0}$ = initial guess; $r_{0} = b - A x_{0}$, $u_0 = A r_0$, $f_0
    = A^T r_0$, $q_0 = v_0 = z_0 = 0$
    \STATE $\sigma_{-1} = \pi_0 = \tau_0 = 0$, $\sigma_0 = (r_0, u_0)$, $\rho_0
    = \alpha_0 = \omega_0 = 1$, $\phi_0 = (r_0, r_0)$
	  \FOR {$j=0,1,\ldots$}
      \STATE $\rho_{j+1} = \phi_j - \omega_j \sigma_{j-1} + \omega_j \alpha_j
      \pi_j$
      \STATE $\delta_{j+1} = \frac{\rho_{j+1}}{\rho_j}\alpha_j$, $\beta_{j+1} =
      \frac{\delta_{j+1}}{\omega_j}$
      \STATE $\tau_{j+1} = \sigma_j + \beta_{j+1} \tau_j - \delta_{j+1}\pi_j$
      \STATE $\alpha_{j+1} = \frac{\rho_{j+1}}{\tau_{j+1}}$
      \STATE $z_{j+1} = \alpha_{j+1} r_j + \beta_{j+1}
      \frac{\alpha_{j+1}}{\alpha_j} z_j - \alpha_{j+1} \delta_{j+1} v_j$, \\
             $v_{j+1} = u_j + \beta_{j+1} v_j - \delta_{j+1} q_j$, \\
             $s_{j+1} = r_j - \alpha_{j+1} v_{j+1}$
      \STATE $q_{j+1} = A v_{j+1}$
      \STATE $t_{j+1} = u_j - \alpha_{j+1} q_{j+1}$, \\
             $\phi_{j+1} = (r_0, s_{j+1})$, \\
             $\pi_{j+1} = (r_0, q_{j+1})$, \\
             $\gamma_{j+1} = (f_0, s_{j+1})$, \\
             $\eta_{j+1} = (f_0, t_{j+1})$, \\
             $\theta_{j+1} = (s_{j+1}, t_{j+1})$, \\
             $\kappa_{j+1} = (t_{j+1}, t_{j+1})$, \\
             $\nu_{j+1} = (s_{j+1}, s_{j+1})$
      \STATE $\omega_{j+1} = \frac{\theta_{j+1}}{\kappa_{j+1}}$
      \STATE $\sigma_{j+1} = \gamma_{j+1} - \omega_{j+1} \eta_{j+1}$
      \STATE $r_{j+1} = s_{j+1} - \omega_{j+1} t_{j+1}$, \\
             $x_{j+1} = x_j + z_{j+1} + \omega_{j+1} s_{j+1}$
      \IF{($\nu_{j+1} - \omega_{j+1} \theta_{j+1}) < \varepsilon^2$}
        \STATE \bf{break}
      \ENDIF
      \STATE $u_{j+1} = A r_{j+1}$
    \ENDFOR
  \end{algorithmic}
\end{algorithm}

%\footnotetext{The listing of Improved BiCGStab algorithm includes two
\footnote[4]{The listing of Improved BiCGStab algorithm includes two
corrections pointed out in the source code of PETSc library~\cite{petsc}.}

\begin{minipage}[t]{0.465\textwidth}
\begin{algorithm}[H]
\small
    \centering
    \caption{The merged Pipelined BiCGStab method}\label{Alg:PipeBiCGStab} 
    \begin{algorithmic}[1]
      \STATE $x_{0}$ = initial guess; $r_{0} = b - A x_{0}$, \\
      $w_0 = A r_0$, $t_0 = A w_0$
      \STATE $\rho_{0} = (r_{0},r_{0})$, $\alpha_0 = \rho_0/(r_0,w_0)$, \\
      $\beta_{-1} = 0$
      
      \FOR {$j=0,1,\ldots$}
        \STATE $p_j = r_j + \beta_{j-1} (p_{j-1} - \omega_{j-1} s_{j-1})$, \\
               $s_j = w_j + \beta_{j-1} (s_{j-1} - \omega_{j-1} z_{j-1})$, \\
               $z_j = t_j + \beta_{j-1} (z_{j-1} - \omega_{j-1} v_{j-1})$, \\
               $q_j = r_j - \alpha_j s_j$, \\
               $y_j = w_j - \alpha_j z_j$, \\
               $\theta_j = (q_j, y_j)$, $\phi_j = (y_j, y_j)$, \\
               $\pi_j = (q_j, q_j)$ 
        \STATE $v_j = A z_j$
        \STATE $\omega_{j} = \frac{\theta_j} {\phi_j}$
        \STATE $x_{j+1} = x_{j} + \alpha_{j} p_j + \omega_j q_j$, \\
               $r_{j+1} = q_j - \omega_j y_j$
        \IF{$(\pi_j - \omega_j \theta_j) < \varepsilon^2$}
          \STATE \bf{break}
        \ENDIF 
        \STATE $w_{j+1} = y_{j} - \omega_{j} (t_j - \alpha_j v_j)$
        \label{PipeBiCGStab_4vec} \\
               $\rho_{j+1} = (r_0, r_{j+1})$, $\psi_j = (r_0, z_{j})$, \\
               $\sigma_j = (r_0, w_{j+1})$, $\delta_j = (r_0, s_{j})$
        \STATE $t_{j+1} = A w_{j+1}$
        \STATE $\beta_{j} = (\alpha_j / \omega_j) (\rho_{j+1} / \rho_j)$
        \STATE $\alpha_{j+1} = \rho_{j+1} / (\sigma_j + \beta_j \delta_j -
        \beta_j \omega_j \psi_j)$
      \ENDFOR
    \end{algorithmic}
\end{algorithm}
\end{minipage}
\begin{minipage}[t]{0.46\textwidth}
\begin{algorithm}[H]
\small
    \centering
    \caption{The merged preconditioned BiCGStab method}
    \begin{algorithmic}[1] 
      \STATE $x_{0}$ = initial guess; $r_{0} = b - A x_{0}$
	  \STATE $\rho_{0} = (r_{0},r_{0})$
	  \STATE $p_{0} = r_{0}$\\
	  \FOR {$j=0,1,\ldots$}
	    \STATE $\hat p_j = M^{-1} p_{j}$, $v_j = A \hat p_j$
	    \STATE $\delta_{j} = (v_j,r_{0})$
	    \STATE $\alpha_{j} = \frac{\rho_{j}}{\delta_j}$
	    \STATE $s_{j} = r_{j} - {\alpha_{j}} v_j$
	    \STATE $\hat s_j = M^{-1} s_{j}$, $t_j = A \hat s_j$
	    \STATE $\phi_j = (t_j, s_{j})$, $\psi_j = (t_j, t_j)$, \\
	           $\theta_j = (s_j, s_j)$
	    \STATE $\omega_{j} = \frac{\phi_j}{\psi_j}$
	    \IF{$(\theta_j - \omega_j \phi_j) < \varepsilon^2$} 
	      \STATE $x_{j+1} = x_{j} + \alpha_{j} \hat p_j + {\omega_{j}} \hat s_j$
	      \STATE $r_{j+1} = s_{j} - {\omega_{j}} t_j$
	      \STATE \bf{break}
	    \ENDIF
	    \STATE $r_{j+1} = s_{j} - {\omega_{j}} t_j$, \\
	  		   $\rho_{j+1} = (r_{j+1},r_{0})$
	  \STATE $\beta_{j} = {\frac{\rho_{j+1}} {\rho_{j}}}
	  {\frac{\alpha_{j}} {\omega_{j}}}$
	  \STATE $x_{j+1} = x_{j} + \alpha_{j} \hat p_{j} + {\omega_{j}} \hat s_{j}$
	  \STATE $p_{j+1} = r_{j+1} + {\beta_{j}} (p_{j}-{\omega_{j}} v_j)$
	  \ENDFOR
    \end{algorithmic}
\end{algorithm}
\end{minipage}
\hfill

\begin{minipage}[t]{0.46\textwidth}
\begin{algorithm}[H]
\small
    \centering
%    \caption{The merged Reordered BiCGStab method\footnotemark} \label{Alg:PRBiCGStab}
    \caption{The merged Reordered BiCGStab method $^5$} \label{Alg:PRBiCGStab}
    \begin{algorithmic}[1] 
      \STATE $x_{0}$ = initial guess; $r_{0} = b - A x_{0}$
      \STATE $\rho_{0} = (r_{0},r_{0})$
      \STATE $z_0 = M^{-1}r_{0}$
      \STATE $\hat v_0 = z_0$
      \FOR {$j=0,1,\ldots$}
        \STATE $v_j = A \hat v_j$
        \STATE $\delta_{j} = (v_j, r_{0})$
        \STATE $s_j = M^{-1} v_j$
        \STATE $\alpha_{j} = \frac{\rho_{j}}{\delta_{j}}$
        \STATE $\hat t_j = z_j - \alpha_{j} s_j$
        \STATE $t_j = A \hat t_j$
        \STATE $\tilde r_j = r_{j} - \alpha_{j} v_j$, \\
               $\theta_j = (t_j, \tilde r_j)$, $\phi_j = (t_j, t_j)$, \\
               $\psi_j = (t_j, r_{0})$, $\eta_j = (\tilde r_j, \tilde r_j)$
         \STATE $q_j = M^{-1} t_j$
        \STATE $\omega_{j} = \frac{\theta_j} {\phi_j}$
        \STATE $r_{j+1} = \tilde r_j - \omega_j t_j$
        \IF{$(\eta_j - \omega_j \theta_j) < \varepsilon^2$}
          \STATE $x_{j+1} = x_{j} + \alpha_{j} \hat v_j + \omega_j \hat t_j$
          \STATE \bf{break}
        \ENDIF 
        \STATE $\rho_{j+1} = - \omega_{j} \psi_j$
        \STATE $\beta_{j} = {\frac{\rho_{j+1}} {\rho_{j}}}
               {\frac{\alpha_{j}} {\omega_{j}}}$
        \STATE $x_{j+1} = x_{j} + \alpha_{j} \hat v_j + \omega_{j} \hat t_j$, \\
               $z_{j+1} = \hat t_j - \omega_{j} q_j$, \label{rbicgstab_typo} \\
               $\hat v_{j+1} = z_{j+1} + \beta_j (\hat v_j -\omega_j s_j)$
      \ENDFOR
    \end{algorithmic}
\end{algorithm}
\end{minipage}
\hfill
\begin{minipage}[t]{0.465\textwidth}
\begin{algorithm}[H]
\small
  \centering
  \caption{The merged preconditioned Pipelined BiCGStab method}
  \label{Alg:PPipeBiCGStab_merged} 
  \begin{algorithmic}[1]
    \STATE $x_{0}$ = initial guess; $r_{0} = b - A x_{0}$, \\
    $\hat r_0 = M^{-1} r_0$, $w_0 = A \hat r_0$, \\
    $\hat w_0 = M^{-1} \hat w_0$, $t_0 = A \hat w_0$ 
    \STATE $\rho_{0} = (r_{0},r_{0})$, $\alpha_0 = \rho_0/(r_0,w_0)$, \\
    $\beta_{-1} = 0$
    \FOR {$j=0,1,\ldots$}
      \STATE $\hat p_j = \hat r_j + \beta_{j-1} (\hat p_{j-1} - \omega_{j-1}
      \hat s_{j-1})$, \\ 
      $\hat s_j = \hat w_j + \beta_{j-1} (\hat s_{j-1} - \omega_{j-1} \hat
      z_{j-1})$, \\
      $\hat q_j = \hat r_j - \alpha_j \hat s_j$ 
      \STATE $s_j = w_j + \beta_{j-1} (s_{j-1} - \omega_{j-1} z_{j-1})$, \\
      $z_j = t_j + \beta_{j-1} (z_{j-1} - \omega_{j-1} v_{j-1})$, \\
      $q_j = r_j - \alpha_j s_j$, \\
      $y_j = w_j - \alpha_j z_j$, \\
      $\theta_j = (q_j, y_j)$, $\phi_j = (y_j, y_j)$, \\
      $\pi_j = (q_j, q_j)$ 
      \STATE $\hat z_j = M^{-1} z_j$, $v_j = A \hat z_j$
      \STATE $\omega_{j} = \frac{\theta_j} {\phi_j}$
      \IF{$(\pi_j - \omega_j \theta_j) < \varepsilon^2$}
        \STATE $x_{j+1} = x_{j} + \alpha_{j} \hat p_j + \omega_j \hat q_j$
        \STATE $r_{j+1} = q_j - \omega_j y_j$
        \STATE \bf{break}
      \ENDIF 
      \STATE $x_{j+1} = x_{j} + \alpha_{j} \hat p_j + \omega_j \hat q_j$, \\
      $\hat r_{j+1} = \hat q_j - \omega_{j} (\hat w_j - \alpha_j \hat z_j)$
      \STATE $r_{j+1} = q_j - \omega_j y_j$, \\
      $w_{j+1} = y_{j} - \omega_{j} (t_j - \alpha_j v_j)$, \\
      $\rho_{j+1} = (r_0, r_{j+1})$, $\psi_j = (r_0, z_{j})$, \\
      $\sigma_j = (r_0, w_{j+1})$, $\delta_j = (r_0, s_{j})$
      \STATE $\hat w_j = M^{-1} w_j$, $t_j = A \hat w_j$
      \STATE $\beta_{j} = (\alpha_j / \omega_j) (\rho_{j+1} / \rho_j)$
      \STATE $\alpha_{j+1} = \rho_{j+1} / (\sigma_j + \beta_j \delta_j -
      \beta_j \omega_j \psi_j)$
    \ENDFOR
  \end{algorithmic}
\end{algorithm}
\end{minipage}
%\footnotetext{The presented listing of the Reordered BiCGStab algorithm fixes a
\footnote[5]{The presented listing of the Reordered BiCGStab algorithm fixes a
typo in the second row of line~\ref{rbicgstab_typo}, introduced in the algorithm
originally presented in~\cite{RBiCGStab}.}

%%%%%%%%%%%%%%%%%%%%%%%%%%%%%%%%%%%%%%%%%%%%%%%%%%%%%%%%%%%

\section {Benchmarking efficiency of asynchronous non-blocking Allreduce
operations}
\label{app:Iallreduce_async}

The efficiency of the asynchronous progression of global communications is among
the key factors influencing on the choice of the optimal iterative method.
The test platform used in the present paper for parallel simulations has
InfiniBand interconnect and does not provide hardware-based progression for
asynchronous execution of non-blocking global communications, thus focusing the
interest on two other progression techniques. The efficiency of the non-blocking
global communications can be measured with help of Intel MPI
benchmark~\cite{IntelMPIBench}, and specifically the IMB-NBC collection of
benchmarks~\cite{Intel_NBC}. These benchmarks allow to compare the execution
times for blocking and non-blocking operations, and estimate the
communication/computation overlap. The basic functionality, however, does not
allow to measure the efficiency of manual progression. To cover all the
simulation scenarios of interest the modified benchmark has been developed. The
proposed benchmark measures three simulation scenarios: blocking global
reduction followed by the computations, non-blocking global reduction overlapped
by the computations, and non-blocking global reduction overlapped by the
computations with manual progression. The corresponding simulation scenarios are
summarized in Figure~\ref{fig:allreduce_bench_proto}. The calculations are
performed in the form of vector updates, and the size of the vectors is chosen
in order to balance the computations with communications. The number of loop
iterations, $N_{it}$, in the experiments varies from 10 to 100 to estimate the
typical number of $MPI\_Test$ calls needed to progress the global communication
as a function of the number of computational processes. The overall amount of
computations is preserved constant by changing the vector size. The experiments
are performed up to 192~compute nodes with 8~MPI processes per each node.

\begin{figure}[ht!]
\begin{lstlisting}
// mode #1:
MPI_Allreduce(...);
for (i = 0; i < N_it; i++)
  vector_op(...);

// mode #2:
MPI_Iallreduce(..., &req);
for (i = 0; i < N_it; i++)
  vector_op(...);
MPI_Wait(&req, ...);

// mode #3:
MPI_Iallreduce(..., &req);
for (i = 0; i < N_it; i++) {
  vector_op(...);
  MPI_Test(&req, ...);
}
MPI_Wait(&req, ...);
\end{lstlisting}
\vspace{-0.5cm}
\caption{Execution scenarios performed in the asynchronous non-blocking
Allreduce communications benchmark.} \label{fig:allreduce_bench_proto}
\end{figure}

\begin{figure}[t!] 
\centering
\includegraphics[width=6cm]{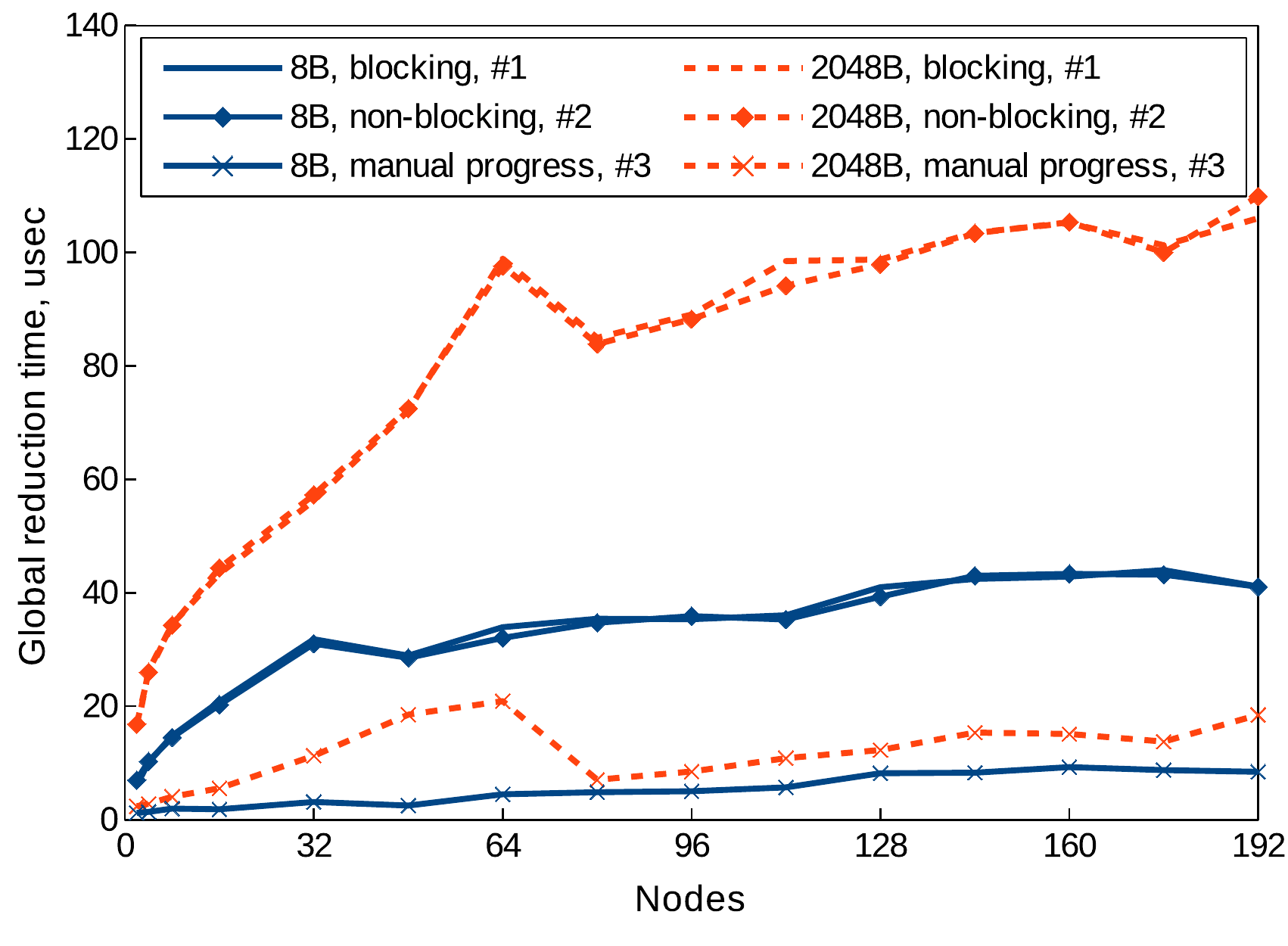}
\includegraphics[width=6cm]{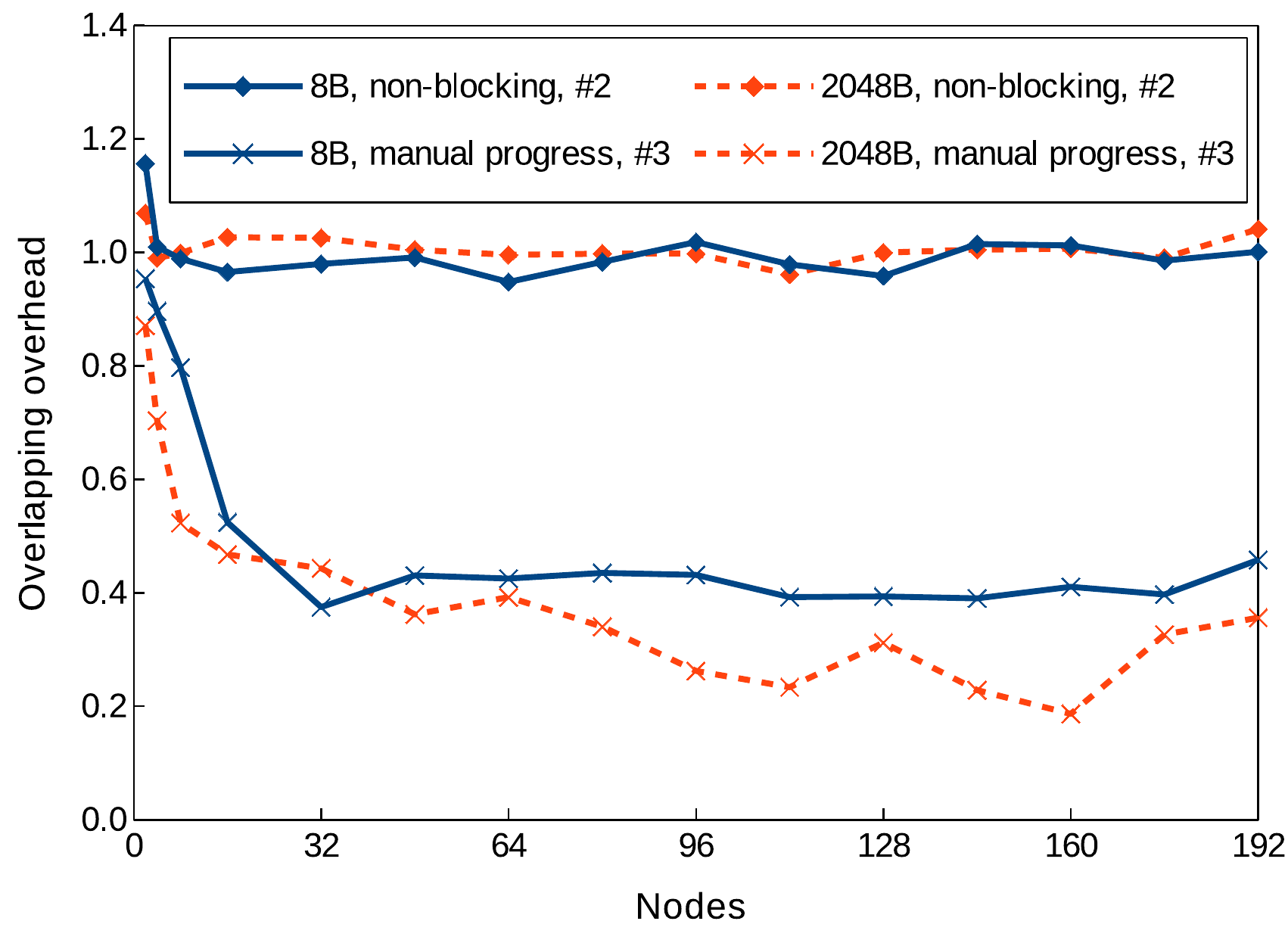}
\vspace{-1.0cm}
\caption{Execution times and overlapping overhead for global reduction
operations with different execution scenarios. Left -- cumulative times for data
transfer functions (\textit{MPI\_Allreduce/MPI\_Iallreduce} and
\textit{MPI\_Wait}); right -- overlapping overhead, $N_{it} = 20$.}
\label{fig:iall_bench}
\end{figure}

The corresponding tests performed on Lomonosov supercomputer with Intel MPI
library 2017 demonstrate that the basic implementation of asynchronous data
transfer without additional progression (mode~\#2) does not allow to obtain any
observable performance gain compared the synchronous operation (mode~\#1).
The execution time for \textit{MPI\_Iallreduce} function is much lower compared
the blocking one, but the total time of \textit{MPI\_Iallreduce} and
\textit{MPI\_Wait} functions is about the same as for \textit{MPI\_Allreduce}
function (Figure~\ref{fig:iall_bench}). Despite the significant time spent on
computations, which theoretically can hide the latency of global communications,
in practice the progress is not automatically performed.

The third execution scenario investigates the effect of \textit{MPI\_Test}
function calls to handle the message progression. While the systematic
\textit{MPI\_Test} calls affect the efficiency of the computations, the manual
progression allows to obtain overlapping of communications and computations. The
corresponding overlapping overhead is estimated by the following parameter
\begin{equation}
\gamma^i = \frac{T_{calc}^{\#i} + T_{comm}^{\#i} -
T_{calc}^{\#1}}{T_{comm}^{\#1}},
\label{eq:allreduce_reduction}
\end{equation}
characterizing ratio of the non-blocking communication time to the blocking one,
and also accounting increase in the calculation time. Results presented in
Figure~\ref{fig:iall_bench} demonstrate that manual progression allows to reduce
the global communications time. Starting from 16~nodes the corresponding
communications time can be reduced at least twofold.

\begin{figure}[t!] 
\centering
\includegraphics[width=6cm]{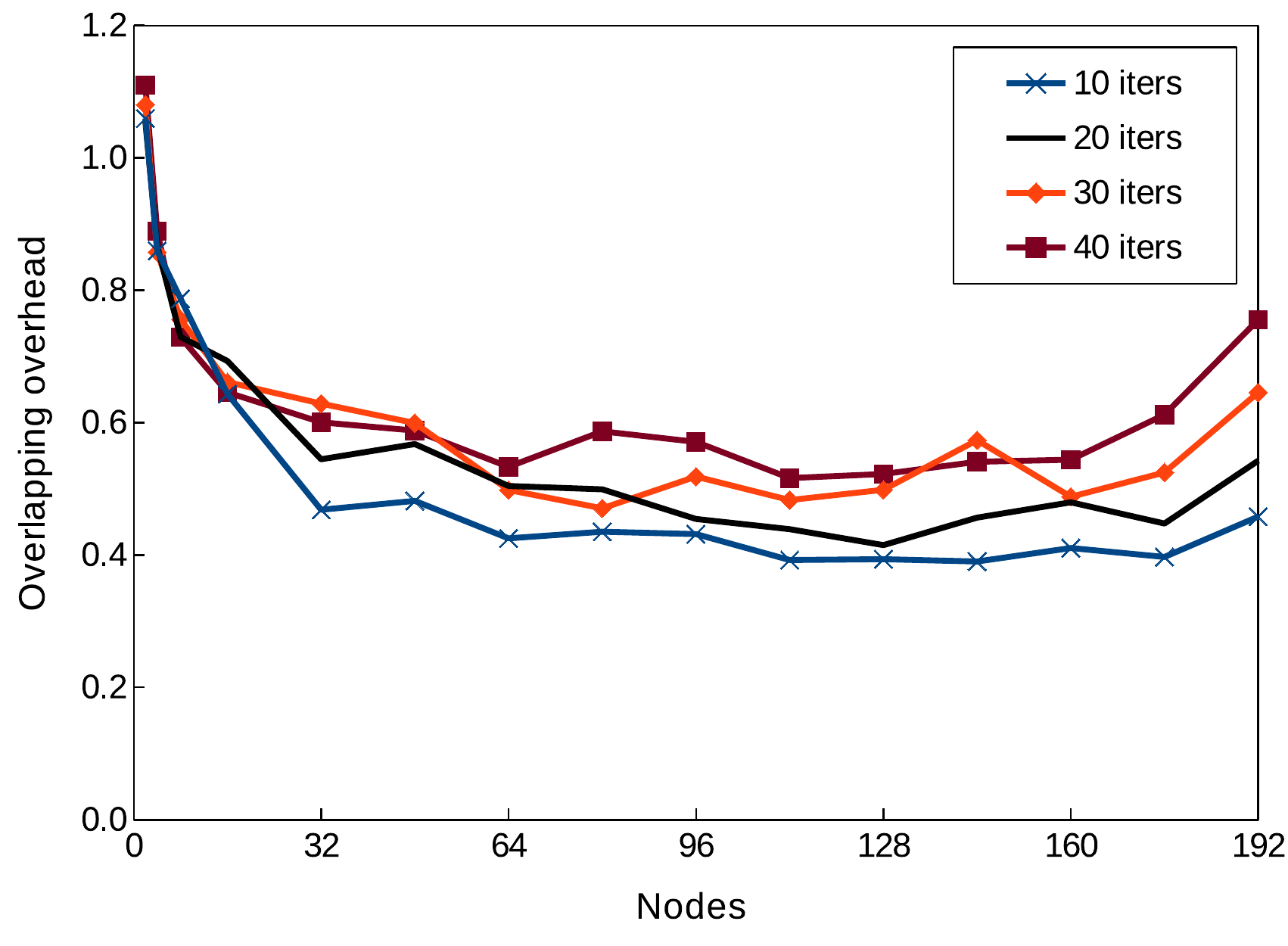}
\includegraphics[width=6cm]{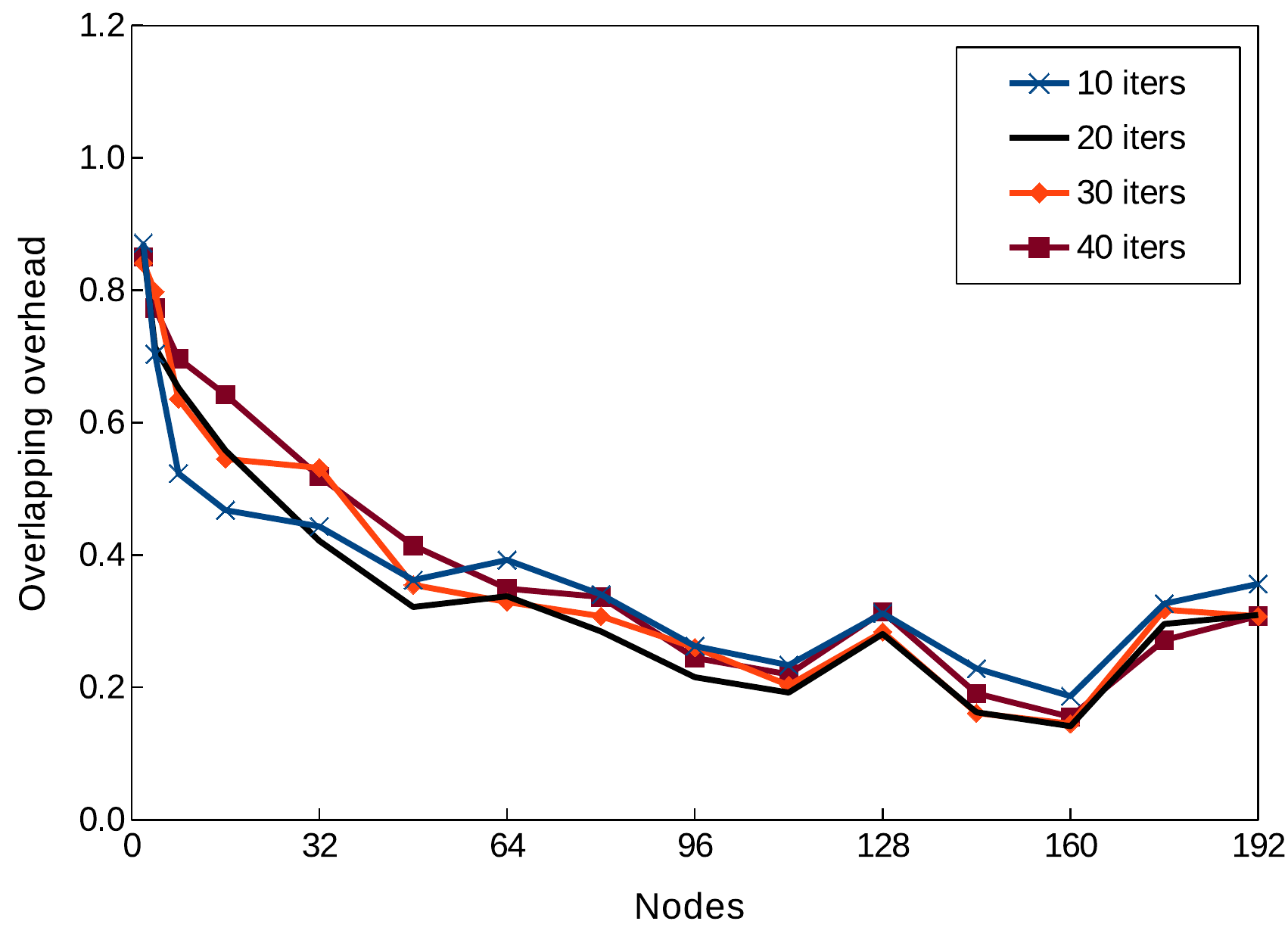}
\vspace{-1.0cm}
\caption{Influence of the number of \textit{MPI\_Test} calls on overlapping
overhead. Left -- message size 8~bytes; right -- message size 2048~bytes.}
\label{fig:iall_ratio}
\end{figure}

The influence of the number of \textit{MPI\_Test} calls on the efficiency of
message progression is shown in Figure~\ref{fig:iall_ratio}. For the short
messages the optimal number of internal loop iterations is close to~10. The
increase in the number of iterations decreases the time spent on final
\textit{MPI\_Wait} function call, but increases the time spent on multiple
\textit{MPI\_Test} calls. The optimal number of iterations grows with the
message size: for the messages of 2048~bytes the corresponding value increases
to 20--30~iterations with overall improvement in the efficiency of
communications overlap.

\begin{figure}[t!]
\centering
\includegraphics[width=6cm]{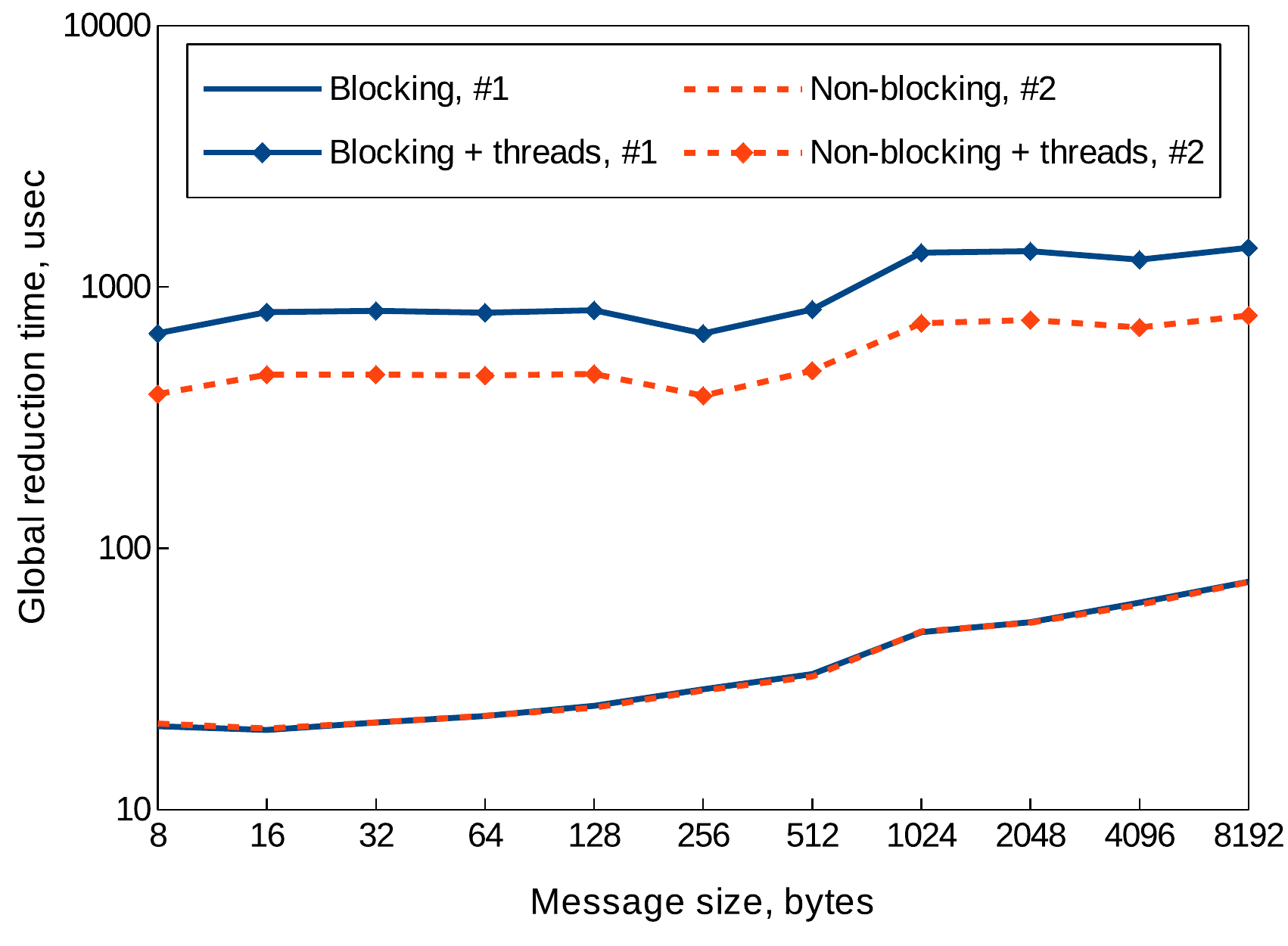}
\includegraphics[width=6cm]{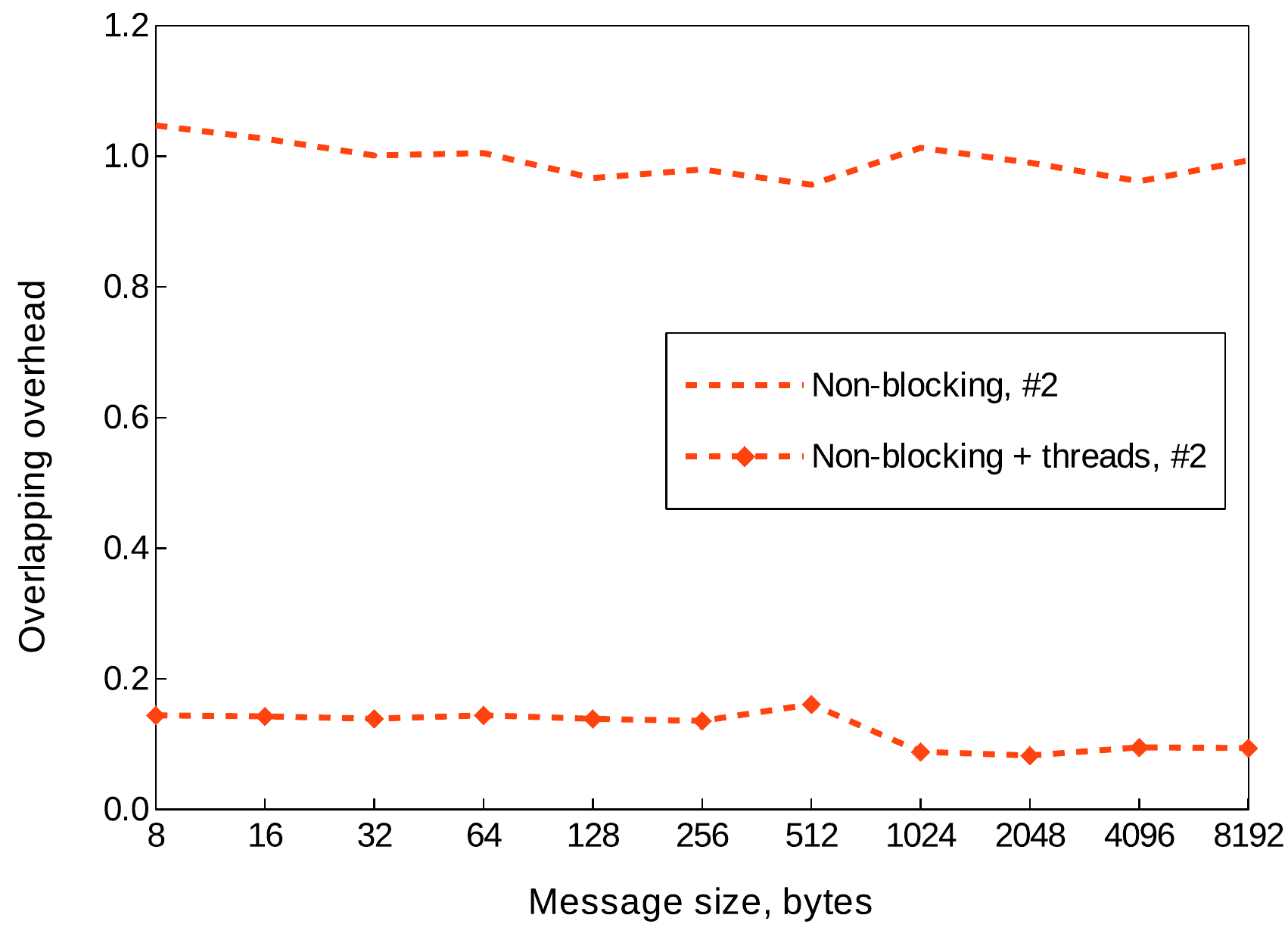}
\vspace{-1cm}
\caption{Execution times and overlapping overhead for global reductions with
asynchronous progression threads. Left -- cumulative times for data transfer
functions (\textit{MPI\_Allreduce/MPI\_Iallreduce} and \textit{MPI\_Wait});
right -- overlapping overhead.}
\label{fig:software_progression}
\end{figure}

Intel MPI library 2017 used in the presented calculations provides software
progression functionality, which can be enabled by setting the environment
variable \textit{I\_MPI\_ASYNC\_PROGRESS = 1}. Activating this option, the
library creates progression threads to handle the asynchronous non-blocking
message progression. The corresponding functionality has been investigated using
the same benchmark described above. Comparison of the obtained results shows
that software progression leads to a significant slowdown of both communications
and computations: the latency of blocking reduction operation increases by a
factor of~40. This observation also correlates with the obtained IMB-NBC
Iallreduce benchmark results and~\cite{IntelMPIBench}. Despite the almost ideal
overlap of communications and computations, the overall time of global reduction
and computations becomes more than by an order of magnitude higher compared the
one without software progression (Figure~\ref{fig:software_progression}).
Reducing the number of MPI processes per node and changing the threads pinning
rule does not allow to obtain any observable improvement in the communication
time. This fact makes the corresponding functionality inapplicable for the
computational algorithms massively performing global reductions with small
vectors.

%%%%%%%%%%%%%%%%%%%%%%%%%%%%%%%%%%%

\section*{Acknowledgements}
The presented work is supported by the RSF grant No. 18-71-10075. 

The author acknowledges Alexey Medvedev for in-depth discussions of materials
presented in the paper. The research is carried out using the equipment of the
shared research facilities of HPC computing resources at Lomonosov Moscow State
University.

%%%%%%%%%%%%%%%%%%%%%%%%%%%%%%%%%%%%%%%%%%%%%%%%%%%%%%%%%%%%%%%%%%%%%%%%%%%%%%%

\bibliography{base_hpc}

%%%%%%%%%%%%%%%%%%%%%%%%%%%%%%%%%%%%%%%%%%%%%%%%%%%%%%%%%%%%%%%%%%%%%%%%%%%%%%%

\end{document}